\documentclass[10pt]{article}

\usepackage{dirtytalk}

\usepackage{amsfonts}
\usepackage{amssymb}
\usepackage{amsmath}
\usepackage{enumerate}
\usepackage{amsthm}

\newcommand{\cf}{\mathrm{cf}}
\newcommand{\cof}{\mathrm{cof}}

\newcommand{\cov}{\mathrm{cov}}

\newcommand{\pp}{\mathrm{pp}}
\newcommand{\PP}{\mathrm{PP}}
\newcommand{\tcf}{\mathrm{tcf}}

\newcommand{\ADS}{\mathrm{ADS}}

\newtheorem{Th}{\bf THEOREM}[section]

\newtheorem{Lem}[Th]{\bf LEMMA}

\newtheorem{Pro}[Th]{\bf PROPOSITION}

\newtheorem{fact}[Th]{\bf FACT}

\newtheorem{Cor}[Th]{\bf COROLLARY}

\newtheorem{Obs}[Th]{\bf OBSERVATION}

\theoremstyle{definition} \newtheorem{Def}[Th]{\bf DEFINITION}
 
\theoremstyle{remark}\newtheorem{Rmk}[Th]{\bf REMARK}

\theoremstyle{question}

\usepackage[english]{babel}
\makeindex

\title{$\mu$-CLUBS OF $P_\kappa (\lambda)$ : PARADISE IN HEAVEN}
\author{Pierre MATET}

\date{}

\begin{document}

\maketitle

\renewcommand{\thefootnote}{\arabic{footnote}} 	

\renewcommand{\thefootnote}{}                                



\vskip 0,7cm

\begin{abstract}  Let $\mu < \kappa < \lambda$ be three infinite cardinals, the first two being regular. We show that if there is no inner model with large cardinals, $u (\kappa, \lambda)$ is regular, where $u (\kappa, \lambda)$ denotes the least size of a cofinal subset in $(P_\kappa (\lambda), \subseteq)$, and $\cf (\lambda) \not= \mu$, then (a) the $\mu$-club filters on $P_\kappa (\lambda)$ and $P_\kappa (u (\kappa, \lambda))$ are isomorphic, and (b) the ideal dual to the $\mu$-club filter on $P_\kappa (\lambda)$ (and hence the restriction of the nonstationary ideal on $P_\kappa (\lambda)$ to sets of uniform cofinality $\mu$) is not $I_{\kappa, \lambda}$-$\frak{b}_{u (\kappa, \lambda)}$-saturated. 
\end{abstract}

\bigskip

\section{Introduction}

\bigskip

\subsection{Sur la terre comme aux cieux}

\bigskip

This is one in a series of papers attempting to show that some of the (near-)perfectness of heaven is reflected on earth. Define \say{earth} : the universe $V$ we live in. Define \say{heaven} : not easy to say, since we are missing some of the details from down below, but this is when the universe is $L$-like in the sense that for instance (a) SSH (Shelah's Strong Hypothesis) holds, and (b) for any singular cardinal $\rho$, weak square (i.e. $\square_\rho^\ast$) holds at $\rho$. So GCH does not necessarily hold, and we see heaven as a federal state, not a centralized one. Roughly speaking, heaven is what the universe is like when there are no inner models with large cardinals (e.g. cardinals $\tau$ of Mitchell order $\tau^{++}$).

The modus operandi is as follows :
\begin{enumerate}[\rm (1)]
\item Start from a statement (A) true in $L$. 
\item Eliminate its component dependent on GCH, thus substituting for (A) a statement (B) that is (a) equivalent to (A) under GCH, and (b) true in heaven.
\item Show that there are situations (for instance for certain values of the cardinals concerned) when (B) holds on earth.
\end{enumerate}

\bigskip

\subsection{Un petit coin de paradis}

\bigskip

There would not be anything wrong with (B) being always true, but we do not expect this to happen too often, since after all earth is not heaven. What we are looking for is just a little corner of paradise.

\bigskip

\subsection{Il n'y a de bon que les choses exquises}

\bigskip

There should be something paradisiac about (B), for instance it could be optimal. 

\bigskip

\subsection{A representative example}

\bigskip

As an example, let us consider the representation problem in pcf theory (missing definitions can be found in Sections 1.6 and 1.7). We are given a regular uncountable cardinal $\kappa$ and a cardinal $\lambda \geq \kappa$. Assuming that $\lambda^{< \kappa} > \lambda$, does there exist a scale of length $\lambda^{< \kappa}$ ? By a result of Shelah (see Fact 1.4 below), the answer is affirmative in $L$, where the following statement holds :

\medskip

(A) Suppose that $\lambda^{< \kappa} > \lambda$. Then $\lambda^{< \kappa} = \tcf (\prod A/I)$ for some infinite set $A$ of regular infinite cardinals with $\vert A \vert  < \kappa < \sup A \leq \lambda$, and some ideal $I$ on $A$ with $\{A \cap a : a \in A\} \subseteq I$.

\medskip

What makes (A) hold in $L$ is that in $L$, if $\lambda^{< \kappa} > \lambda$, then $\cf (\lambda) < \kappa$ and $\lambda^{< \kappa} = \lambda^+$. (A) is dependent on cardinal exponentiation. In fact it will fail in a generic extension of $L$ obtained by adding, say, $\lambda^{++}$ Cohen reals. To obtain a GCH-free version of (A), we replace $\lambda^{< \kappa}$ with its pcf theory \say{core}, namely $u (\kappa, \lambda)$ (the least size of a cofinal subset of $P_\kappa (\lambda)$) :

\medskip

(B) Suppose that $u (\kappa, \lambda) > \lambda$. Then $u (\kappa, \lambda) = \tcf (\prod A/I)$ for some infinite set $A$ of regular infinite cardinals with $\vert A \vert  < \kappa < \sup A \leq \lambda$, and some ideal $I$ on $A$ with $\{A \cap a : a \in A\} \subseteq I$.

\medskip

(B) is both equivalent to (A) under GCH and true under SSH, but so is also

\medskip

(B') Suppose that $u (\kappa, \lambda) > \lambda$. Then $\lambda^+ = \tcf (\prod A/I)$ for some infinite set $A$ of regular infinite cardinals with $\vert A \vert  < \kappa < \sup A \leq \lambda$, and some ideal $I$ on $A$ with $\{A \cap a : a \in A\} \subseteq I$.

\medskip

But (B'), being much weaker than (B), does not have the taste of paradise, so we will stick to (B). By results of Shelah (see Fact 1.14 below), (B) will hold if the following conditions are satisfied :
\begin{enumerate}[\rm (a)]
\item $\kappa$ and $u (\kappa, \lambda)$ are both successor cardinals.
\item $\tau^{< \kappa} < \lambda$ for every cardinal $\tau < \lambda$.
\item Either $\cf (\lambda) \not= \omega$, or $\lambda$ is not a fixed point of the aleph function.
\end{enumerate}
With condition (b), cardinal exponentiation (which we were trying to keep at bay) is back ! More in the spirit of pcf theory, it is shown in \cite{pcf} that (B) holds if (a) and (b') both hold, where (b') asserts that $\lambda$ is smaller than the least point of the aleph function above $\kappa$. How to go further ? One thing we could do would be to replace in (B) the condition that $\sup A \leq \lambda$ with the weaker condition that $\sup A < u (\kappa, \lambda)$. But the main problem is elsewhere : $u (\kappa, \lambda)$ might be just too big. What makes us suspicious is that if $A$ and $I$ are as in the statement of (B), then by a result of Shelah, we must have that 

\centerline{$u (\kappa, \lambda) \leq \pp_{\vert A \vert} (\sup A) \leq \cov(\sup A, \sup A, \vert A \vert^+, 2)$.} 

So in \cite{Menas} we considered the variant (B") obtained by replacing in (B) $u (\kappa, \lambda)$ with a new quantity denoted by $\pp (\kappa, \lambda)$ and defined by

\centerline{$\pp (\kappa, \lambda) = \max (\lambda, \sup \{\pp (\chi) : \cf (\chi) < \kappa < \chi \leq \lambda \})$.}

Obviously, (B") will hold provided that $\pp (\kappa, \lambda)$ is a successor cardinal. It is shown in \cite{Menas} that $\pp (\kappa, \lambda)$ is less than or equal to $u (\kappa, \lambda)$, and equal to it in case $\kappa$ is a successor cardinal and $\lambda$ is less than the least point of the aleph function above $\kappa$. 

\medskip

We see this example as a good illustration of the difficulties involved in finding the right formulation for version (B) of an original statement (A) true in $L$.

\bigskip

\subsection{The secret life of $\mu$-clubs}

\bigskip

Let $\mu < \kappa < \lambda$ be three infinite cardinals, the first two being regular. There are quite a few ideals on $P_\kappa (\lambda)$ that can be seen as generalizations to $P_\kappa (\lambda)$ of $NS_\kappa \vert E^\kappa_\mu$ (the restriction of the nonstationary ideal on $\kappa$ to the set of ordinals less than $\kappa$ of cofinality $\mu$). Let us mention some of them (for missing definitions see Section 2) :
\begin{itemize}
\item $NS_{\kappa, \lambda} \vert E^{\kappa, \lambda}_\mu$ (the restriction of the nonstationary ideal on $P_\kappa (\lambda)$ to the set of elements of $P_\kappa (\lambda)$ of uniform cofinality $\mu$.)
\item $NS_{\mu, \kappa, \lambda} \vert E^{\kappa, \lambda}_\mu$, where $NS_{\mu, \kappa, \lambda}$ denotes the smallest $(\mu, \kappa)$-normal ideal on $P_\kappa (\lambda)$.
\item $N\mu$-$S_{\kappa,\lambda}$ (the ideal dual to the filter on $P_\kappa (\lambda)$ generated by the $\mu$-club subsets of $P_\kappa (\lambda)$, where $C \subseteq P_\kappa (\lambda)$ is a $\mu$-club if it is cofinal and closed under unions of $\subset$-increasing sequences of length $\mu$).
\item $NG^\mu_{\kappa, \lambda}$ (the definition of which involves games on $P_\kappa (\lambda)$ of length $\mu$).
\end{itemize}

Which ones of these ideals are mere restrictions of $NS_{\kappa, \lambda}$ ? More generally, how do these ideals compare ? We addressed these questions in \cite{Resemble} and \cite{Secret}. In particular we showed there that it may happen that $N\mu$-$S_{\kappa,\lambda}$ is so large that it cannot be a restriction of  $NS_{\kappa, \lambda}$. Largeness of an ideal can be measured by the size of its set of generators. It can also be witnessed by an isomorphism between the ideal and an ideal on a bigger set. We will establish the following : 

\medskip

\begin{Th}  Suppose that there is no inner model with large cardinals, $u (\kappa, \lambda)$ is regular, and $\cf (\lambda) \not= \mu$. Then $N\mu$-$S_{\kappa,\lambda}$ and $N\mu$-$S_{\kappa, u (\kappa, \lambda)}$ are isomorphic.
\end{Th}

\medskip

Notice the so-called critical cofinality $\cf (\lambda) = \mu$.

\medskip

Now if $N\mu$-$S_{\kappa,\lambda}$ is isomorphic to $N\mu$-$S_{\kappa, u (\kappa, \lambda)}$, it inherits the properties of $N\mu$-$S_{\kappa, u (\kappa, \lambda)}$, in particular in terms of non-saturation. We will study non-saturation properties of $N\mu$-$S_{\kappa, \sigma}$ for $\sigma (= (u (\kappa, \lambda))$ regular, and then apply Theorem 1.1, which will give us the following :

\medskip

\begin{Th}  Suppose that there is no inner model with large cardinals, $u (\kappa, \lambda)$ is regular, and $\cf (\lambda) \not= \mu$. Then there exists an $(N\mu$-$S_{\kappa,\lambda}, I_{\kappa, \lambda})$-tower of length $\frak{b}_{u (\kappa, \lambda)}$ (and hence $N\mu$-$S_{\kappa,\lambda}$ is not $I_{\kappa, \lambda}$-$\frak{b}_{u (\kappa, \lambda)}$-saturated). 
\end{Th}

\bigskip

\subsection{Heaven}

\bigskip

Theorems 1.1 and 1.2 both assume that there is no inner model with large cardinals. Let us be more specific. We need to recall some notions of Shelah's pcf theory (see \cite{SheCA}, \cite{HSW}).

\medskip

\begin{Def} Let $A$ be an infinite set of regular infinite cardinals such that $\vert A \vert  < \min (A)$, and $I$ be an ideal on $A$. For $f, g \in {}^A On$, let $f <_I g$ if $\{ a \in A : f(a) < g(a)\} \in I^\ast$,  and  $f\leq_I g$ if $\{ a\in A : f(a) \leq g(a)\} \in I^\ast$.

For $h \in {}^A  On$, let $\prod h =   \displaystyle\prod_{a\in A} h(a)$.  We let $\prod A = \prod  h$,  where $h(a) = a$ for all $a\in A$.

Let  $h : A \rightarrow On \setminus \{ 0\}$,  and $\pi$ be a regular cardinal. We say that $\tcf (\prod h/I) = \pi$ if there exists an $<_I$-increasing sequence $\vec{f} = \langle f_\xi : \xi < \pi\rangle$ in $\prod h$ such that for any $g \in \prod h$, there is $\xi < \pi$ with $g \leq_I f_\xi$.  Such a sequence is said to be a {\it scale of length} $\pi$ for $\sup (A)$.

For a singular cardinal $\sigma$, we let 

\centerline{$\pp (\sigma) = \sup (\{ \pi : \pi$ is the length of some scale for $\sigma \})$.}
\end{Def}

\begin{fact}  {\rm (\cite[Theorem 1.5 p. 50]{SheCA})} Let $\sigma$ be a singular cardinal. Then there is a set $A$ of regular cardinals such that o.t.$(A) = \cf(\sigma) < \min A$, $\sup A = \sigma$ and $\tcf (\prod A / I) = \sigma^+$, where $I$ is the noncofinal ideal on $A$.
\end{fact}

\begin{Def}  {\it Shelah's Strong Hypothesis} {\rm (SSH)} asserts that $\pp(\sigma) = \sigma^+$ for every singular cardinal $\sigma$.
\end{Def}

\medskip

The exact consistency strength of the failure of SSH is not known, but it is thought to be roughly the same as that of the failure of SCH (the Singular Cardinal Hypothesis) that has been shown by Gitik \cite{Gitik}) to be equiconsistent with the existence of a measurable cardinal $\chi$ of Mitchell order $\chi^{++}$.

Let us next recall Shelah's Almost Disjoint Set principle.

\medskip

\begin{Def} Let $\tau \leq \theta$ be two infinite cardinals. $\ADS_\tau (\theta)$ asserts the existence of a sequence $\langle z_\beta : \beta < \theta \rangle$  such that 
\begin{itemize}
\item for any $\beta < \theta$, $z_\beta$ is a cofinal subset of $\tau$ of order type  $\cf(\tau)$ ;
\item for any $e \in P_{\tau^+} (\theta)$, there is $g : e \rightarrow \tau$ with the property that $(z_\beta \setminus g(\beta)) \cap (z_\gamma \setminus g(\gamma)) = \emptyset$ whenever $\beta < \gamma$ are in $e$.
\end{itemize}
\end{Def}

\begin{fact} {\rm(\cite{CFM}, see also \cite{Scales})} Let $\sigma$ be a singular cardinal such that $\square_\sigma^\ast$ holds. Then $\ADS_\sigma (\sigma^+)$ holds.
\end{fact}

\medskip

Thus the consistency strength of the failure of $\ADS_\sigma (\sigma^+)$ at some singular cardinal $\sigma$ is quite high (the failure of $\square_\sigma^\ast$ for $\sigma$ singular is known to entail a Woodin cardinal (see \cite[p. 702]{Jech2}). 

Our assumption for Theorems 1.1 and 1.2 is that SSH holds, and moreover $\ADS_\sigma (\sigma^+)$ holds for every singular cardinal $\sigma$. The assumption will hold if 
we are in heaven (in the sense of Section 1.1). Everything in there is perfect, so there is not much left to do except describing this very perfection and endlessly rejoicing about being part of it (which may explain why  Eve found the place fairly boring after a while.)

\medskip

In a companion paper \cite{OnEarth}, we will show that if $u (\kappa, \lambda)$ is regular, $\cf (\lambda) \not= \mu$, and $\kappa$ and $\lambda$ assume certain values, then $N\mu$-$S_{\kappa,\lambda}$ and $N\mu$-$S_{\kappa, u (\kappa, \lambda)}$ are isomorphic, and moreover there exists an $(N\mu$-$S_{\kappa,\lambda}, I_{\kappa, \lambda})$-tower of length $\frak{b}_{u (\kappa, \lambda)}$.

\bigskip

\subsection{GCH stripped of its clothes}

\bigskip

In a way, SSH constitutes the hard core of GCH. To see this, let us first recall Shelah's covering numbers.

\medskip

\begin{Def}  Given four cardinals $\rho_1, \rho_2, \rho_3, \rho_4$ with $\rho_1 \geq \rho_2 \geq \rho_3 \geq \omega$ and $\rho_3 \geq \rho_4 \geq 2$, $\cov (\rho_1, \rho_2, \rho_3, \rho_4)$ denotes the least cardinality of any $X \subseteq  P_{\rho_2}(\rho_1)$ such that for any $a \in  P_{\rho_3}(\rho_1)$, there is $Q \in  P_{\rho_4}(X)$ with $a \subseteq \bigcup Q$.          
\end{Def}

\begin{fact}   {\rm (\cite[pp. 85-86]{SheCA}, \cite{LCCN})}  Let $\rho_1, \rho_2, \rho_3$ and $\rho_4$ be four cardinals such that $\rho_1 \geq \rho_2 \geq \rho_3 \geq \omega$ and $\rho_3 \geq \rho_4 \geq 2$. Then the following hold :        
\begin{enumerate}[\rm (i)]
\item  If $\rho_1 = \rho_2$ and either $\cf(\rho_1) < \rho_4$ or $\cf(\rho_1) \geq \rho_3$, then $\cov (\rho_1, \rho_2, \rho_3, \rho_4) = \cf(\rho_1)$.           
\item  If either $\rho_1 > \rho_2$, or $\rho_1 = \rho_2$ and $\rho_4 \leq \cf(\rho_1) < \rho_3$, then $\cov (\rho_1, \rho_2, \rho_3, \rho_4) \geq \rho_1$.        
\item  $\cov (\rho_1, \rho_2, \rho_3, \rho_4) = \cov (\rho_1, \rho_2, \rho_3, \max \{\omega, \rho_4\})$.
\item $\cov (\rho_1^+, \rho_2, \rho_3, \rho_4) = \max \{\rho_1^+, \cov (\rho_1, \rho_2, \rho_3, \rho_4)\}$.
\item  If $\rho_1 > \rho_2$ and $\cf(\rho_1) < \rho_4 = \cf(\rho_4)$, then 

\centerline{$\cov (\rho_1, \rho_2, \rho_3, \rho_4) = \sup \{\cov (\rho, \rho_2, \rho_3, \rho_4) : \rho_2 \leq \rho < \rho_1\}$.}

\item   If $\rho_1$ is a limit cardinal such that $\rho_1 > \rho_2$ and $\cf(\rho_1) \geq \rho_3$, then 

\centerline{$\cov (\rho_1, \rho_2, \rho_3, \rho_4) = \sup \{\cov (\rho, \rho_2, \rho_3, \rho_4) : \rho_2 \leq \rho < \rho_1\}$.}

\item   If $\rho_3 > \rho_4 \geq \omega$, then
 
\centerline{$\cov (\rho_1, \rho_2, \rho_3, \rho_4) = \sup \{\cov (\rho_1, \rho_2, \rho^+, \rho_4) : \rho_4 \leq \rho < \rho_3\}$.}
       
\item   If $\rho_3 \leq \rho_2 = \cf(\rho_2)$, $\omega \leq \rho_4 = \cf(\rho_4)$ and $\rho_1 < \rho_2^{+\rho_4}$, then $\cov (\rho_1, \rho_2, \rho_3, \rho_4) = \rho_1$. 
\item   If $\rho_3 = \cf(\rho_3)$, then either  $\cf(\cov (\rho_1, \rho_2, \rho_3, \rho_4)) < \rho_4$, or  $\cf(\cov (\rho_1, \rho_2, \rho_3, \rho_4)) \geq \rho_3$.
 \end{enumerate}
\end{fact}

\medskip

Under SSH, $\cov (\rho_1, \rho_2, \rho_3, \rho_4)$ can be computed in all cases.

\medskip

\begin{fact} {\rm (\cite{Secret})} Let $\rho_1, \rho_2, \rho_3$ and $\rho_4$ be four infinite cardinals such that $\rho_1 \geq \rho_2 \geq \rho_3 \geq \rho_4$. Assume the SSH. Then the following hold :        
\begin{enumerate}[\rm (i)]
\item  If $\rho_1 = \rho_2$ and either $\cf(\rho_1) < \rho_4$ or $\cf(\rho_1) \geq \rho_3$, then $\cov (\rho_1, \rho_2, \rho_3, \rho_4) = \cf(\rho_1)$.           
\item  If $\rho_4 \leq \cf (\rho_1) < \rho_3$, then $\cov (\rho_1, \rho_2, \rho_3, \rho_4) = \rho_1^+$.        
\item  In all other cases, $\cov (\rho_1, \rho_2, \rho_3, \rho_4) = \rho_1$.
 \end{enumerate}
\end{fact}


\begin{Def}  Given two uncountable cardinals $\tau = \cf (\tau) \leq \sigma$, we let $u (\tau, \sigma) = \cov (\sigma, \tau, \tau, 2)$.          
\end{Def}

\begin{Rmk} $u (\tau, \sigma)$ is a pcf-theoretic version of $\sigma^{< \tau}$, obtained, as often in pcf theory, by replacing $=$ with $\subseteq$ (since $u (\tau, \sigma) =$ the cofinality of the poset $(P_\tau (\sigma), \subseteq)$, whereas $\sigma^{< \tau} =$ the cofinality of the poset $(P_\tau (\sigma), =)$). It is not the only one (see e.g. $\pp (\tau, \sigma)$ defined in \cite{Menas}).
\end{Rmk}

\medskip

Shelah showed that SSH can be reformulated in terms of covering numbers.

\medskip

\begin{fact} {\rm (\cite[p. 59]{SheCA}, \cite{She93}, see also \cite{LCCN})} The following are equivalent :
 \begin{enumerate}[\rm (i)]
 \item  SSH holds.       
\item  $\pp(\sigma) = \sigma^+$ for every singular cardinal $\sigma$ of cofinality $\omega$.
 \item Given two uncountable cardinals $\chi$ and $\nu$ such that $\cf(\chi) = \chi \leq \nu$, $u (\chi, \nu)$ equals $\nu$ if $\cf (\nu) \geq \chi$, and $\nu^+$ otherwise. 
\end{enumerate}
\end{fact}






\begin{fact} \begin{enumerate}[\rm (i)] 
\item {\rm (See e.g. \cite[Theorems 9.1.2 and 9.1.3]{HSW})} Let $\sigma$ be a singular cardinal such that $\rho^{\cf (\sigma)} < \sigma$ for any cardinal $\rho < \sigma$. Suppose that either $\cf (\sigma) \not= \omega$, or $\sigma$ is not a fixed point of the aleph function. Then $\pp (\sigma) = \sigma^{\cf (\sigma)}$.  
\item  {\rm (\cite[Conclusion 1.8 p. 369]{SheCA})} Let $\sigma$ be a singular cardinal of cofinality $\omega$. Suppose that
\begin{itemize}
\item $\pp (\tau) < \sigma$ for any singular cardinal $\tau < \sigma$ of cofinality $\omega$.
\item $\pp (\nu) = \nu^+$ for any large enough singular cardinal $\nu < \sigma$ of cofinality $\omega_1$.
\end{itemize}
Then $\pp (\sigma) = \cov (\sigma, \sigma, \omega_1, 2)$.
\end{enumerate}
\end{fact} 

\begin{Obs} The following are equivalent :
\begin{enumerate}[\rm (i)] 
\item GCH holds.
\item  SSH holds, and so does GCH at regular cardinals (i.e. $2^\tau = \tau^+$ for any regular cardinal $\tau$).
\end{enumerate}
\end{Obs}

{\bf Proof.}

 \hskip0,4cm  (i) $\rightarrow$ (ii) : Trivial.
  
\hskip0,2cm  (ii) $\rightarrow$ (i) : Suppose that (ii) holds, and let $\sigma$ be a singular cardinal. Then clearly, $\sigma$ is a strong limit cardinal. If either $\cf (\sigma) \not= \omega,$ or $\sigma$ is not a fixed point of the aleph function, then by Fact 1.14 (i), $2^\sigma = \sigma^{\cf (\sigma)} = \pp (\sigma) = \sigma^+$. Now suppose that $\sigma$ is a fixed point of the aleph function of cofinality $\omega$. Then by Fact 1.14 (ii), $ \cov (\sigma, \sigma, \omega_1, 2) = \pp (\sigma) = \sigma^+$. Since $\rho^{\aleph_0} \leq 2^\rho < \sigma$ for any infinite cardinal $\rho < \sigma$, it easily follows that $2^\sigma = \sigma^{\aleph_0} \leq \sigma^+$. 
 \hfill$\square$  

\medskip

Of course, GCH at regular cardinals is something that can be easily destroyed by (Easton) forcing \cite{East}. Thus SSH can be seen as a naked body, with the power function on regular cardinals depending on the clothes used to dress it.

\bigskip

\subsection{Surviving without GCH}

\bigskip

To illustrate our point, let us discuss the following result proved below (Corollary 3.34).

\medskip

\begin{Pro} Suppose that GCH holds, $\cf (\lambda) < \mu$, and $\kappa$ is the successor of a cardinal $\nu$ with $\cf (\nu) \not= \cf (\lambda)$. Then $N\mu$-$S_{\kappa,\lambda}$ and $N\mu$-$S_{\kappa, \lambda^{< \kappa}}$ are isomorphic. 
\end{Pro}

\medskip

A shadow version of this remains true (at least in some cases) if we replace GCH with SSH (Proposition 3.33).

\medskip

\begin{Pro} Suppose that SSH holds, $\cf (\lambda) < \mu$, and $\kappa$ is the successor of a cardinal $\nu$ with $\cf (\nu) \not= \cf (\lambda)$. Suppose further that $(\cf (\lambda))^{++} < \nu$. Then $N\mu$-$S_{\kappa,\lambda}$ and $N\mu$-$S_{\kappa, u (\kappa, \lambda)}$ are isomorphic. 
\end{Pro}

\medskip

Notice that we had to reformulate the conclusion, substituting $u (\kappa, \lambda)$ for $\lambda^{< \kappa}$, this move being justified by the fact that $u (\kappa, \lambda)$ and $\lambda^{< \kappa}$ are equal under GCH.

\medskip

The missing cases can be dealt with by throwing in an additional fragment of GCH (Propositions 3.7 and 3.27).

\medskip

\begin{Pro} Suppose that SSH holds, $\cf (\lambda) < \mu$, and $\kappa$ is the successor of a cardinal $\nu$ with $\cf (\nu) \not= \cf (\lambda)$. Suppose further that $P_\kappa (\lambda)$ carries a pseudo-Kurepa family of size $u (\kappa, \lambda)$. Then $N\mu$-$S_{\kappa,\lambda}$ and $N\mu$-$S_{\kappa, u (\kappa, \lambda)}$ are isomorphic. 
\end{Pro}

\medskip

We think that there should be other GCH results of set theory that could be processed in this fashion.

\bigskip

\subsection{Meanwhile on earth}

\bigskip

Our appeal to SSH was made a posteriori, to make our results look optimal. The original results read as follows.

\medskip

\begin{Pro}  Suppose that $\lambda$ is regular. Then there exists an $(N\mu$-$S_{\kappa,\lambda}, I_{\kappa, \lambda})$-tower of length $\frak{b}_\lambda$. 
\end{Pro}

{\bf Proof.} By Propositions 4.24 and 4.32.
\hfill$\square$

\begin{Pro}  Suppose that $\cf (\lambda) \in \kappa \setminus \{ \mu\}$, and one of the following conditions is satisfied : 
\begin{enumerate}[\rm (a)]
\item $\kappa$ is weakly inaccessible.
\item $\kappa$ is a successor cardinal, and letting $\rho$ denote the largest limit cardinal below $\kappa$, either $\cf (\rho) \not= \mu$, or $\cf (\rho) = \mu < \rho$ and $\pp (\rho) < \kappa$, or $(\max (\rho, \mu))^{+ 3} < \kappa$, or $\rho^\mu < \kappa$.
\item There is no inner model with large large cardinals. 
\end{enumerate}
Then the following hold :
\begin{enumerate}[\rm (i)]
\item $N\mu$-$S_{\kappa,\lambda}$ and $N\mu$-$S_{\kappa, \lambda^+}$ are isomorphic. 
\item There exists an $(N\mu$-$S_{\kappa,\lambda}, I_{\kappa, \lambda})$-tower of length $\frak{b}_{\lambda^+}$. 
\end{enumerate}
\end{Pro}

{\bf Proof.} By Fact 3.31 and Propositions 3.52 and 4.34 (i).
\hfill$\square$


\bigskip

\subsection{Organization}

\bigskip

The paper is organized as follows. In Section 2 we recall definition and properties of the various ideals on $P_\kappa (\lambda)$ that we will be studying. Section 3 is devoted to the isomorphism problem. In Subsection 3.1 we introduce shuttles. They allow us to travel between $P_\kappa (\lambda)$ and $P_\kappa (u (\kappa, \lambda))$. Subsections 3.2 and 3.3 review material concerning pseudo-Kurepa families which are used in Subsections 3.4 and 3.6 to construct shuttles. In Subsection 3.5 we show that in the easier case $\cf (\lambda) < \mu$, we can get not only $N\mu$-$S_{\kappa,\lambda}$ and $N\mu$-$S_{\kappa, u (\kappa, \lambda)}$ to be isomorphic, but also $NS_{\mu, \kappa,\lambda}$ and $NS_{\kappa, u (\kappa, \lambda)}$, or $NG^\mu_{\kappa,\lambda}$ and $NG^\mu_{\kappa,u (\kappa, \lambda)}$ . The results of Section 3 are summarized in Subsection 3.7. Section 4 is concerned with the non-saturation problem. Some preparation is needed, namely results on the sup function on $P_\kappa (\lambda)$ in Subsection 4.1, and results on towers in Subsection 4.2. The main result is stated in Subsection 4.6. Its proof is split in three cases which are respectively dealt with in Subsections 4.3 - 4.5. Two cases are left, when $\kappa \leq \cf (\lambda) < \lambda$ and when $\cf (\lambda) = \mu$, for which we obtain partial results in Section 5.

\bigskip

\section{Ideals}

\bigskip

In this section we review basic definitions and facts concerning the ideals we will be working with.

\medskip

\begin{Def} For a set $A$ and a cardinal $\chi$, we set $P_\chi (A) = \{ a\subseteq A :  \vert a \vert < \chi\}$ and $ [A]^\chi = \{ a\subseteq A :  \vert a \vert  = \chi\}$.
\end{Def}

\begin{Def} Let $X$ be an infinite set. An {\it ideal  on} $X$ is a nonempty collection $J$ of subsets of $X$ such that 
 \begin{itemize}
\item $X \notin J$.
\item $P(A)\subseteq J$ for all $A\in J$.  
\item $A\cup B\in J$ whenever $A, B \in J$. 
\end{itemize}

Given an ideal $J$ on $X$, we let $J^+ = P(X) \setminus J$, $J^\ast = \{ A\subseteq X : X \setminus A\in J\}$, and $J \vert A = \{ B\subseteq X : B \cap A\in J\}$  for each $A\in J^+$.  
For a cardinal   $\rho,J$ is {\it $\rho$-complete} if  $\bigcup Q \in J$ for every $Q \subseteq J$ with $\vert Q \vert< \rho$. ${\rm cof} (J)$ denotes the least cardinality of any $Q \subseteq J$ such that $J = \bigcup_{A\in Q} P(A)$. For a cardinal $\rho$ and $Y\subseteq P(X)$, $J$ is {\it $Y$-$\rho$-saturated} if there is no $Q \subseteq J^+$ with $\vert Q \vert = \rho$ such that $A \cap B \in Y$ for any two distinct members $A, B$ of $Q$. 
We say that $J$ is {\it $\rho$-saturated} 
if it is $J$-$\rho$-saturated 
.  
\end{Def}

 \begin{Def} Let $X$ and $Y$ be two infinite sets, and $J$ be an ideal  on $X$. Given $f : X \rightarrow Y$,  we let $f (J) = \{ B \subseteq Y : f^{-1} (B) \in J\}$.   
\end{Def} 

\begin{fact} {\rm (\cite{Secret})} Suppose that $f^{-1} (\emptyset) \in J$. Then 
\begin{enumerate}[\rm (i)]
\item $f (J)$ is an ideal on $Y$.
\item If $J$ is $\rho$-complete, then so is $f (J)$.
\item $(f (J))^\ast = \{C \subseteq Y : f^{-1} (C) \in J^\ast\}$.
\item $\{ f`` D : D \in J^\ast\} \subseteq (f (J))^\ast$.
\item $\cof (f (J)) \leq \cof (J)$.
\end{enumerate}
\end{fact} 

\begin{Def}  For $i = 1, 2$, let $X_i$  be an infinite set, and $K_i$ be an ideal  on $X_i$. We say that $K_1$ is {\it isomorphic to} $K_2$ if there are $W_1 \in K_1^\ast$, $W_2 \in K_2^\ast$ and a bijection $k : W_1 \rightarrow W_2$ such that
\begin{itemize}
\item for any $A \in K_1^\ast \cap P (W_1)$, $k``A \in K_2^\ast$ ; 
\item for any $B \in K_2^\ast \cap P (W_2)$, there is a unique $A \in K_1^\ast \cap P (W_1)$ such that $B = k``A$.
\end{itemize}
\end{Def}

\begin{Rmk} Beware that for some authors, $K_1$ is isomorphic to $K_2$ if there is $k : W_1 \rightarrow W_2$ as above for $W_1 = X_1$ and $W_2 = X_2$.
\end{Rmk}

\begin{fact} {\rm (Folklore)} Suppose that $K_1$ is isomorphic to $K_2$. Then $K_2$ is isomorphic to $K_1$.
\end{fact} 

\begin{Obs} Suppose that $K_1$ is isomorphic to $K_2$. Then $\cof (K_1) = \cof (K_2)$.
\end{Obs} 

{\bf Proof.} By Fact 2.7, it suffices to establish that $\cof (K_1) \geq \cof (K_2)$. Let $W_1$, $W_2$ and $k$ be as in Definition 2.5. Notice that $k^{-1} (\emptyset) = \emptyset$. For $i = 1, 2$, put $J_i = K_i \cap P (W_i)$. Clearly, $J_i$ is an ideal on $W_i$. Moreover, $\cof (J_i) \leq \cof (K_i) \leq 1 + \cof (J_i)$. Since the cofinality of an ideal is either 1, or infinite, it follows that $\cof (J_i) = \cof (K_i)$.

 \medskip

{\bf Claim.}   $k (J_1) = J_2$.      

\medskip

{\bf Proof of the claim.} 

$\supseteq$ : Let $B \in J_2$. Since $W_2 \setminus B \in K_2^\ast$, there is $A \in K_1^\ast \cap P (W_1)$ such that $B = k``A$. Then $A = k^{- 1} (W_2 \setminus B)$, so

\centerline{$W_1 \setminus k^{- 1} (W_2 \setminus B) = k^{- 1} (B)$}

lies in $K_1$ and hence in $J_1$.

$\subseteq$ : Let $X \in k (J_1)$. Since $W_1 \setminus k^{- 1} (X) \in K_1^\ast$, we must have that $k`` (W_1 \setminus k^{- 1} (X)) \in K_2^\ast$. Then

\centerline{$W_2 \setminus k`` (W_1 \setminus k^{- 1} (X)) = W_2 \setminus (W_2 \setminus X) = X$}

lies in $K_2$ and hence in $J_2$, which completes the proof of the claim. 

\medskip 

By the claim and Fact 2.4, 

\centerline{$\cof (K_2) = \cof (J_2) \leq \cof (J_1) = \cof (K_1)$.}

\hfill$\square$ 

\begin{Obs} The following are equivalent :
\begin{enumerate}[\rm (i)]
\item  $f (K_1) = K_2$ for some one-to-one $f : X_1 \rightarrow X_2$.
\item  There are $W_2 \in K_2^\ast$ and a bijection $k : X_1 \rightarrow W_2$ such that
\begin{itemize}
\item for any $A \in K_1^\ast$, $k``A \in K_2^\ast$ ; 
\item for any $B \in K_2^\ast \cap P (W_2)$, there is a unique $A \in K_1^\ast$ such that $B = k``A$.
\end{itemize}
\end{enumerate}
\end{Obs} 

{\bf Proof.} 

 \hskip0,4cm  (i) $\rightarrow$ (ii) : Given a one-to-one $f : X_1 \rightarrow X_2$ such that $K_2 = f (K_1)$, put $W_2 = ran (f)$, and define $k : X_1 \rightarrow W_2$ by $k (a) = f (a)$. Now for any $A \in K_1^\ast$, we have that $A \subseteq f^{- 1} (f``A)$, so $f^{- 1} (f``A) \in K_1^\ast$. It follows that $f``A (= k``A)$ lies  in $K_2^\ast$. Finally, let $B \in K_2^\ast \cap P (W_2)$. Put $G = h^{- 1} (B)$. Then obviously, $G \in K_1^\ast$, and moreover $k``G = f``G = B$. Suppose toward a contradiction that there is $H \in K_1^\ast$ such that $H \not= G$ and $B = k``H$. Pick $a \in H \bigtriangleup K$. Then since $k$ is one-to-one, $k (a) \in k``H \bigtriangleup k``G$. Contradiction ! 
 
 \medskip

\hskip0,2cm  (ii) $\rightarrow$ (i) : Suppose that $W_2$ and $k$ are in (ii). 

\medskip

{\bf Claim 1.}  $K_2 \subseteq k (K_1)$.      

\medskip

{\bf Proof of Claim 1.} Let $C \in K_2^\ast$. Clearly, $C \cap W_2 \in K_2^\ast$, so we may find $A \in K_1^\ast$ such that $C \cap W_2 = k``A$. Then $A \subseteq k^{- 1} (C)$, and consequently $k^{- 1} (C) \in K_1^\ast$, which completes the proof of the claim. 

\medskip 

{\bf Claim 2.}  $k (K_1) \subseteq K_2$.      

\medskip

{\bf Proof of Claim 2.} Let $C \subseteq X_2$ be such that $k^{- 1} (C) \in K_1^\ast$.  Then $k``(k^{- 1} (C)) \in K_2^\ast$. Since $k``(k^{- 1} (C)) \subseteq C$, it follows that $C \in K_2^\ast$, which completes the proof of the claim. 

\hfill$\square$ 

\medskip

Thus if $f (K_1) = K_2$ for some one-to-one $f : X_1 \rightarrow X_2$, then $K_1$ and $K_2$ are isomorphic.

\medskip

\begin{Def}  For a regular infinite cardinal $\chi$, $I_\chi$ and $NS_\chi$ denote, respectively, the noncofinal ideal on $\chi$ and  the nonstationary ideal on $\chi$.

Given a regular uncountable cardinal $\chi$ and a cardinal $\sigma \geq \chi$, $I_{\chi, \sigma}$ and $NS_{\chi, \sigma}$ denote, respectively,  the noncofinal ideal on $P_\chi (\sigma)$ and the nonstationary ideal on $P_\chi (\sigma)$.  
 
 \medskip
 
 An ideal $J$ on $P_\chi (\sigma)$ is {\it fine} if $I_{\chi,\sigma} \subseteq J$. 
\end{Def}

\begin{fact} {\rm (Folklore)} Let $\tau > \sigma$ be a cardinal. Then $I_{\chi, \sigma} = q (I_{\chi, \tau})$,  where $q : P_\chi (\tau) \rightarrow P_\chi (\sigma)$ is defined by $q (x) = x \cap \sigma$.
\end{fact}



\medskip

Let us recall some properties of $NS_{\chi, \sigma}$.

\medskip

\begin{Def}  For $F : \sigma \times \sigma \rightarrow \sigma$, we let $C_F$ denote the set of all nonempty $a \in P_\chi (\sigma)$ such that $a \cap \chi \in \chi$ and $F``(a \times a) \subseteq a$.
\end{Def}  

\medskip

The following is well-known (see e.g. \cite{Conc}).

\begin{fact} Let $C \subseteq P_\chi (\sigma)$. Then the  following are equivalent :
 \begin{enumerate}[\rm (i)]
 \item  $C \in NS_{\chi, \sigma}^\ast$.   
 \item  There is $F : \sigma \times \sigma \rightarrow \sigma$ such that $C_F \subseteq C$.   
 \end{enumerate}
\end{fact}

\begin{fact} {\rm (\cite{Men})} Let $\tau > \sigma$ be a cardinal. Then $NS_{\chi, \sigma} = q (NS_{\chi, \tau})$,  where $q : P_\chi (\tau) \rightarrow P_\chi (\sigma)$ is defined by $q (x) = x \cap \sigma$.
\end{fact}


\medskip Slightly more basic than the nonstationary ideal is the ideal dual to the strong club filter.

\medskip

\begin{Def}  A subset $C$ of $P_\chi (\sigma)$ is {\it strongly closed} if $\bigcup X \in C$ for all $X \in P_\chi (C) \setminus \{\emptyset\}$. 

$SNS_{\chi, \sigma}$ denotes the collection of all $B \subseteq P_\chi (\sigma)$ such that $B \cap C = \emptyset$ for some strongly closed, cofinal subset $C$ of $P_\chi (\sigma)$.
\end{Def} 

\medskip

\begin{fact} {\rm (\cite{Conc})} Suppose $\sigma > \chi$. Then the following hold :
 
\begin{enumerate}[\rm (i)]
\item  $SNS_{\chi, \sigma}$ is a $\chi$-complete, fine ideal on $P_\chi (\sigma)$.
\item  $I_{\chi, \sigma} \subset SNS_{\chi, \sigma} \subset NS_{\chi, \sigma}$. 
\item Let $j$ be a  one-to-one function from $\sigma \times \sigma$ to $\sigma$, and $D$ be the set of all $a \in P_\chi (\sigma)$ such that $j (\alpha, \beta) \in a$ whenever $\alpha < \beta$ are in $a$. Then $NS_{\chi, \sigma} = SNS_{\chi, \sigma} \vert D$.
\end{enumerate}
\end{fact}


\begin{fact} {\rm (\cite{Secret})} $SNS_{\chi, \sigma} \vert X = NS_{\chi, \sigma} \vert X$ for any $X \in SNS^+_{\chi, \sigma}$ such that $NS_{\chi, \sigma} \subseteq SNS_{\chi, \sigma} \vert X$.
\end{fact}

\begin{Def}  For $f : \sigma \rightarrow P_\chi (\sigma)$, let 

\centerline{$D_f = \{a \in P_\chi (\sigma) \setminus \{\emptyset\} : f``a \subseteq P(a)\}$.}
\end{Def}  

\medskip

\begin{fact} {\rm (\cite{Men})} Let $C \subseteq P_\chi (\sigma)$. Then the  following are equivalent :
 \begin{enumerate}[\rm (i)]
 \item  $C \in SNS_{\chi, \sigma}^\ast$.   
 \item  There is $f : \sigma \rightarrow P_\chi (\sigma)$ such that $D_f \subseteq C$.   
 \end{enumerate}
\end{fact}

\begin{fact} {\rm (\cite{Secret})} Let $\theta$ be a regular cardinal with $\kappa \leq \theta \leq \sigma$, and $C$ be a closed unbounded subset of $\theta$. Then $\{a \in P_\kappa (\sigma) : \sup (a \cap \theta) \in C\} \in SNS_{\kappa, \sigma}^\ast$.
\end{fact}

\medskip

Let us next turn to a stronger notion of normality. $(\nu, \chi)$-normality is a light version of the property of $[\chi]^{< \nu}$-normality studied in \cite{MPS1}. The advantage of $(\nu, \chi)$-normal ideals is that they always exist.

\medskip

\begin{Def} Let $\nu$ be a regular cardinal less than $\chi$. A $\chi$-complete, fine ideal $K$ on $P_\chi (\sigma)$ is $(\nu, \chi)${\it -normal} if for any $A \in K^+$ and any $f : A \rightarrow P_\nu (\sigma)$ with the property that $f (a) \subseteq a$ for every $a \in A$, there are $d \in P_\chi (\sigma)$ and $B \in K^+ \cap P (A)$ such that $f (a) \subseteq d$ for all $a \in B$.
\end{Def}
 
\begin{fact} {\rm(\cite{Resemble}, \cite{Secret})}
\begin{enumerate}[\rm (i)]
\item If $K$ is $(\nu, \chi)$-normal, then it is normal.
\item Assume $\nu = \omega$. Then $K$ is $(\nu, \chi)$-normal if and only if it is normal.
\item There exists a smallest  $(\mu, \chi)$-normal ideal on $P_\chi (\sigma)$.
\item  Let $J$ be a $(\nu, \chi)$-normal ideal on $P_\chi (\tau)$, where $\tau$ is a cardinal greater than $\sigma$. Then $q (K)$ is $(\nu, \chi)$-normal, where $q : P_\chi (\tau) \rightarrow P_\chi (\sigma)$ is defined by $q (x) = x \cap \sigma$.        
\end{enumerate}
\end{fact}

\begin{Def} We let $NS_{\nu, \chi, \sigma}$ denote the smallest $(\nu, \chi)$-normal ideal on $P_\kappa (\sigma)$.
\end{Def}

\medskip

Some of our results will involve  game ideals.

\medskip

\begin{Def} For $A \subseteq P_\chi (\sigma)$, $G_{\chi,\sigma}^\nu (A)$ denotes the following two-person game consisting of $\nu$ moves. At step  $\alpha < \nu$, player I selects $a_\alpha \in P_\chi (\sigma)$, and II replies by playing $b_\alpha \in P_\chi (\sigma)$. The players must follow the rule that for $\beta < \alpha < \nu$, $b_\beta \subseteq a_\alpha \subseteq b_\alpha$. II wins if and only if $\bigcup_{\alpha < \nu} a_\alpha \in A$.  

$NG_{\chi,\sigma}^\mu$ denotes the collection of all $A \subseteq P_\chi (\sigma)$ such that II has a winning strategy in $G_{\chi,\sigma}^\nu (P_\chi (\sigma) \setminus A)$.
 \end{Def}

\begin{fact} 
\begin{enumerate}[(i)]
\item  {\rm (\cite{Secret})} $NG_{\chi,\sigma}^\nu$ is a $(\nu, \chi)$-normal ideal on $P_\chi (\sigma)$.
\item  {\rm (\cite{Ideals})} Let $\tau > \sigma$ be a cardinal. Then $NG_{\chi, \sigma}^\nu =  q (NG_{\chi,\tau}^\nu)$,  where $q : P_\chi (\tau) \rightarrow P_\chi (\sigma)$ is defined by $q (x) = x \cap \sigma$.
\end{enumerate}
\end{fact}

\medskip

Our main object of study is the (ideal dual to) the $\nu$-club filter.

\medskip

\begin{Def} A subset $A$ of $P_\chi (\sigma)$ is $\nu${\it -closed} if $\bigcup_{i < \nu} a_i \in A$ for every increasing sequence $\langle a_i : i < \nu \rangle$ in $(A, \subset)$.   

A subset $C$ of $P_\chi (\sigma)$ is a $\nu${\it -club} if it is a $\nu$-closed, cofinal subset of $P_\chi (\sigma)$. 

We let  $N\nu$-$S_{\chi,\sigma}$ be the set of all $B \subseteq P_\chi(\sigma)$ such that $B \cap C = \emptyset$ for some $\nu$-club $C\subseteq P_\chi (\sigma)$. 
\end{Def}

\begin{fact} \begin{enumerate}[\rm (i)]
\item  {\rm(\cite{Conc})}  $N\nu$-$S_{\chi,\sigma}$ is a normal, fine ideal on $P_\chi (\sigma)$.
\item   {\rm(\cite{Ideals})} $N\nu$-$S_{\chi,\sigma} \subseteq NG_{\chi,\sigma}^\nu$.
\item {\rm (\cite{Secret})}  Let $\tau > \sigma$ be a cardinal. Then $N\nu$-$S_{\chi,\sigma} \subseteq q (N\nu$-$S_{\chi,\tau})$,  where $q : P_\chi (\tau) \rightarrow P_\chi (\sigma)$ is defined by $q (x) = x \cap \sigma$
.\end{enumerate}
\end{fact} 

\medskip

We conclude our tour with the restriction of the nonstationary ideal to the set of all $a$ of uniform cofinality $\nu$.

\medskip

\begin{Def} $E_\nu^\chi$ (respectively, $E_{\geq \nu}^\chi$) denotes the set of all limit ordinals $\alpha < \chi$ with $\cf (\alpha) = \nu$ (respectively, $\cf (\alpha) \geq \nu$).

We let $E_\nu^{\chi, \sigma}$ denote the set of all $a \in P_\chi (\sigma)$ such that for any cardinal $\rho \leq \sigma$ of cofinality greater than or equal to $\chi$, $\sup(a \cap \rho)$ is a limit ordinal of cofinality $\nu$ that does not belong to $a$. 
\end{Def}  


\begin{fact} {\rm(\cite{Ideals})} $NS_{\chi, \sigma} \vert E_\nu^{\chi, \sigma} \subseteq N\nu$-$S_{\chi,\sigma}$.
\end{fact}

\bigskip

\section{Isomorphisms} 

\bigskip

In this section we look for isomorphisms between $N\mu$-$S_{\kappa,\lambda}$ and $N\mu$-$S_{\kappa,\pi}$ for a cardinal $\pi \geq \lambda$. To travel between $P_\kappa (\lambda)$ and $P_\kappa (\pi)$, we will use shuttles.

\medskip

\subsection{Shuttles}

\bigskip

\begin{Def}  Let $\pi \geq \lambda$ be a cardinal. A $(\mu, \kappa, \lambda, \pi)$-{\it shuttle} is a function $\chi : P_\kappa (\lambda) \rightarrow P_\kappa (\pi)$ such that
\begin{enumerate}[\rm (a)]
\item $\chi (a) \cap \lambda = a$ for all $a \in P_\kappa (\lambda)$.
\item $\chi (a) \subseteq \chi (b)$ whenever $a, b \in P_\kappa (\lambda)$ are such that $a \subseteq b$.
\item $ran (\chi) \in I^+_{\kappa, \pi}$.
\item $\chi (\bigcup_{i < \mu} a_i) \subseteq \bigcup_{i < \mu} \chi (a_i)$ for any increasing sequence $\langle a_i : i < \mu \rangle$ in $(P_\kappa (\lambda), \subset)$.
\end{enumerate}
\end{Def}

\begin{Rmk} In case $\pi = \lambda$, the identity function on $P_\kappa (\lambda)$ is a $(\mu, \kappa, \lambda, \pi)$-shuttle, and it is the only one.
\end{Rmk}

\begin{Obs} Let $\chi : P_\kappa (\lambda) \rightarrow P_\kappa (\pi)$ be a $(\mu, \kappa, \lambda, \pi)$-shuttle. Then the following hold :
\begin{enumerate}[\rm (i)]
\item  $\chi$ is one-to-one.
\item $a \subseteq b$ whenever $a, b \in P_\kappa (\lambda)$ are such that $\chi (a) \subseteq \chi (b)$.
\item $\chi (a) \subset \chi (b)$ whenever $a, b \in P_\kappa (\lambda)$ are such that $a \subset b$.
\item $\chi`` A \in I^+_{\kappa, \pi}$ for all $A \in I^+_{\kappa, \lambda}$.
\item $\chi (\bigcup_{i < \mu} a_i) = \bigcup_{i < \mu} \chi (a_i)$ for any increasing sequence $\langle a_i : i < \mu \rangle$ in $(P_\kappa (\lambda), \subset)$.
\end{enumerate}
\end{Obs}

\begin{Cor} Suppose that there exists a $(\mu, \kappa, \lambda, \pi)$-shuttle. Then $u (\kappa, \lambda) = u (\kappa, \pi)$.
\end{Cor}

\medskip

By Fact 1.9, it follows that if there exists a $(\mu, \kappa, \lambda, \pi)$-shuttle, then $\pi \leq u (\kappa, \lambda)$.

\medskip

\begin{Obs} Let $\chi : P_\kappa (\lambda) \rightarrow P_\kappa (\pi)$ be a $(\mu, \kappa, \lambda, \pi)$-shuttle. Then $\chi (I_{\kappa, \lambda}) = I_{\kappa, \pi} \vert ran (\chi)$.
\end{Obs}

{\bf Proof.} It is simple to see that ($I_{\kappa, \pi}$ and hence) $I_{\kappa, \pi} \vert ran (\chi)$ is included in $\chi (I_{\kappa, \lambda})$. For the reverse inclusion, fix $X \in \chi (I_{\kappa, \lambda})$. There must be $b \in P_\kappa (\lambda)$ such that $b \setminus c \not= \emptyset$ for any $c \in P_\kappa (\lambda)$ with $\chi (c) \in X$. Now let $x \in X \cap ran (X)$. Put $x = \chi (c)$. Then since $c = \chi (c) \cap \lambda = x \cap \lambda$, we have that ($b \setminus (x \cap \lambda)$ and hence) $b \setminus x$ is nonempty.
\hfill$\square$

\begin{Lem}  Let $\chi : P_\kappa (\lambda) \rightarrow P_\kappa (\pi)$ be a $(\mu, \kappa, \lambda, \pi)$-shuttle, and $D$ be a $\mu$-club subset of $P_\kappa (\pi)$. Then $\chi^{- 1} (D) \in I^+_{\kappa, \lambda}$.
\end{Lem}

{\bf Proof.} Given $a \in P_\kappa (\lambda)$, inductively define $x_i, w_i \in D$ for $i < \mu$ so that
\begin{itemize}
\item $\chi (a) \subseteq x_0 = w_0$.
\item $x_i = w_i = \bigcup_{j < i} w_j$ in case $i$ is an infinite limit ordinal.
\item $x_i \subseteq w_{i + 1} \subset x_{i + 1}$.
\item $\chi (x_i \cap \lambda) \subseteq w_{i + 1} \subseteq \chi (x_{i + 1} \cap \lambda)$.
\end{itemize}
Finally, put $x = \bigcup_{i < \mu} x_i$. Then the following are readily checked :
\begin{itemize}
\item $x \in D$.
\item $a \subseteq x \cap \lambda$.
\item $\chi (x \cap \lambda) = \chi (\bigcup_{i < \mu} (x_i \cap \lambda)) = \bigcup_{i < \mu} \chi (x_i \cap \lambda) = \bigcup_{i < \mu} w_i = x$.
\end{itemize}
Hence $x \cap \lambda \in \chi^{- 1} (D)$.
\hfill$\square$ 

\begin{Pro} Let $\chi : P_\kappa (\lambda) \rightarrow P_\kappa (\pi)$ be a $(\mu, \kappa, \lambda, \pi)$-shuttle. Then $\chi (N\mu$-$S_{\kappa,\lambda}) = N\mu$-$S_{\kappa,\pi}$ (and hence by Fact 2.7 and Observation 2.9, $N\mu$-$S_{\kappa,\lambda}$ and $N\mu$-$S_{\kappa,\pi}$ are isomorphic).
\end{Pro}

{\bf Proof.} $\supseteq$ : Given $C \in (N\mu$-$S_{\kappa,\pi})^\ast$, pick a $\mu$-club subset $D$ of $C$. By Lemma 3.6, $\chi^{- 1} (D) \in I^+_{\kappa, \lambda}$. To show that $\chi^{- 1} (D)$ is $\mu$-closed, fix $a_i \in \chi^{- 1} (D)$ for $i < \mu$ such that $a_i \subset a_j$ whenever $i < j < \mu$. Then clearly, $\chi (a_i) \subset \chi (a_j)$ whenever $i < j < \mu$, so $\bigcup_{i < \mu} \chi (a_i) \in D$. Since $\bigcup_{i < \mu} \chi (a_i) = \chi (\bigcup_{i < \mu} a_i)$, it follows that $\bigcup_{i < \mu} a_i \in \chi^{- 1} (D)$. Thus $\chi^{- 1} (C) \in (N\mu$-$S_{\kappa,\lambda})^\ast$.

\medskip

$\subseteq$ : Given $C \subseteq P_\kappa (\pi)$ such that $\chi^{- 1} (C) \in (N\mu$-$S_{\kappa,\lambda})^\ast$, select a $\mu$-club subset $D$ of $\chi^{- 1} (C)$. Then $\chi`` D \subseteq C$, and moreover by Observation 3.3 (iv), $\chi`` D \in I^+_{\kappa, \pi}$. To prove that $\chi`` D$ is $\mu$-closed, fix $a_i \in D$ for $i < \mu$ such that $\chi (a_i) \subset \chi (a_j)$ whenever $i < j < \mu$. Then by Observation 3.3 (ii), $a_i \subset a_j$ whenever $i < j < \mu$, so $\bigcup_{i < \mu} a_i \in D$. It follows that $\chi (\bigcup_{i < \mu} a_i)$ (and hence $\bigcup_{i < \mu} \chi (a_i)$) lies in $\chi`` D$. Thus $C \in (N\mu$-$S_{\kappa,\pi})^\ast$.
\hfill$\square$

\begin{Def}  We define $p_\pi : P_\kappa (\pi) \rightarrow P_\kappa (\lambda)$ by $p_\pi (x) = x \cap \lambda$.
\end{Def} 

\medskip

We recalled earlier (Facts 2.11, 2.14 and 2.25 (ii)) that $p_\pi (I_{\kappa,\pi}) = I_{\kappa,\lambda}$, $p_\pi (NS_{\kappa,\pi}) = NS_{\kappa,\lambda}$ and $p_\pi (NG^\mu_{\kappa,\pi}) = NG^\mu_{\kappa,\lambda}$. The following provides another result of this type (with the difference that this one is conditional).

\medskip

\begin{Cor} Let $\chi : P_\kappa (\lambda) \rightarrow P_\kappa (\pi)$ be a $(\mu, \kappa, \lambda, \pi)$-shuttle. Then $p_\pi (N\mu$-$S_{\kappa,\pi}) = N\mu$-$S_{\kappa,\lambda}$.
\end{Cor}

{\bf Proof.} $\supseteq$ : By Fact 2.27 (iii).

$\subseteq$ : Let $C \subseteq P_\kappa (\lambda)$ be such that $p_\pi^{- 1} (C) \in (N\mu$-$S_{\kappa,\pi})^\ast$. Set $D = \chi^{- 1} (p_\pi^{- 1} (C))$. Then by Proposition 3.7, $D \in (N\mu$-$S_{\kappa,\lambda})^\ast$. It is simple to see that $D \subseteq C$.
\hfill$\square$

\bigskip

\subsection{Pseudo-Kurepa families}

\bigskip

Shuttles will be obtained from pseudo-Kurepa families. Pseudo-Kurepa and Kurepa families in $P_\kappa (\lambda)$ with $\kappa = \omega_1$ were introduced by Todorcevic in \cite{Todor}. We find it more convenient to work with sequences. We need some definitions.

\medskip

\begin{Def} Let $\tau$ be a cardinal with $2\leq\tau\leq\kappa$. A $(\tau,\lambda,\pi)${\it -sequence} is a one-to-one sequence $\vec y = \langle y_\beta : \beta < \pi \rangle$ of elements of $P_\tau (\lambda)$  with the property that $y_\beta = \{\beta\}$  for every $\beta < \lambda$.

Given a $(\tau,\lambda,\pi)$-sequence $\vec y$, define $f_{\vec y } : P_\kappa (\lambda) \rightarrow P (\pi)$ by $f_{\vec y } (a) = \{\beta < \pi : y_\beta \subseteq a\}$, and let

\centerline{$A_\kappa (\vec y) = \{x \in P_\kappa (\pi) : f_{\vec y } (x \cap \lambda) \subseteq x\}$}

and 
 
\centerline{$\Delta_\kappa (\vec y) = \{x \in P_\kappa (\pi) : \forall \delta \in x (y_\delta \subseteq x)\}$.} 

\end{Def}

\begin{fact} {\rm (\cite{Secret})}
\begin{enumerate}[\rm (i)]
\item $a =  f_{\vec y} (a) \cap \lambda$.
\item $A_\kappa (\vec y) \cap \Delta_\kappa (\vec y) = \{x \in P_\kappa (\pi) : f_{\vec y} (x \cap \lambda) = x\} = ran ( f_{\vec y})$.
\item $\Delta_\kappa (\vec y) \in SNS^\ast_{\kappa, \pi}$. 
\end{enumerate}
\end{fact}

\begin{Obs} The following are equivalent :
\begin{enumerate}[\rm (i)]
\item $A_\kappa (\vec y) \in I^+_{\kappa, \pi}$.
\item $\vert f_{\vec y} (a) \vert < \kappa$ for every $a \in P_\kappa (\lambda)$. 
\end{enumerate}
\end{Obs}

{\bf Proof.} 

(i) $\rightarrow$ (ii) : Assume that $A_\kappa (\vec y) \in I^+_{\kappa, \pi}$. Given $a \in P_\kappa (\lambda)$, pick $t \in A_\kappa (\vec y)$ with $a \subseteq t$. Then $f_{\vec y} (a)  \subseteq f_{\vec y} (t \cap \lambda)  \subseteq t$, so $\vert f_{\vec y} (a) \vert \leq \vert t \vert< \kappa$.

\medskip

(ii) $\rightarrow$ (i) : Suppose that $\vert f_{\vec y} (a) \vert < \kappa$ for all $a \in P_\kappa (\lambda)$. Given $x \in P_\kappa (\pi)$, set $t = x \cup  f_{\vec y} (x \cap \lambda) = x \cup (f_{\vec y} (x \cap \lambda) \setminus \lambda)$. Notice that $\vert t \vert < \kappa$. Clearly, any $\beta \in f_{\vec y} (t \cap \lambda)$ lies in $f_{\vec y} (x \cap \lambda)$ and hence in $t$. Thus $t \in A_\kappa (\vec y)$.
\hfill$\square$

\begin{Def} An ${\mathcal A}_{\tau,\lambda} (\kappa,\pi)${\it -sequence} is a a $(\tau,\lambda,\pi)$-sequence $\vec y$ with the property that $\vert f_{\vec y} (a) \vert < \kappa$ for every $a \in P_\kappa (\lambda)$.

${\mathcal A}_{\tau,\lambda} (\kappa,\pi)$ asserts the existence of an ${\mathcal A}_{\tau,\lambda} (\kappa,\pi)$-sequence.
\end{Def}

\begin{Obs} Let $\vec y = \langle y_\delta : \delta < \pi \rangle$ be an ${\mathcal A}_{\kappa,\lambda} (\kappa,\pi)$-sequence. Then the following hold : 
\begin{enumerate}[\rm (i)]
\item $ran ( f_{\vec y}) \in I_{\kappa, \pi}^+$.
\item $f_{\vec y} $ is an isomorphism from $(P_\kappa (\lambda), \subset)$ to $(ran ( f_{\vec y}), \subset)$.
\item $I_{\kappa, \pi} \vert ran ( f_{\vec y}) = f_{\vec y} (I_{\kappa, \lambda})$ (and hence by Fact 2.7 and Observation 2.9, $I_{\kappa, \lambda}$ and $I_{\kappa, \pi} \vert ran (f_{\vec y})$ are isomorphic).
\end{enumerate}
\end{Obs}

{\bf Proof.} (i) : Given $t \in P_\kappa (\pi)$, set $b = \bigcup_{\beta \in t} y_\beta$. Then $t \subseteq f_{\vec y} (b)$.

\medskip

(ii) : Use Fact 3.11 (a).

\medskip

(iii) : Fix $X \subseteq P_\kappa (\pi)$.

\medskip

{\bf Claim 1.}  Suppose that $X \notin f_{\vec y} (I_{\kappa, \lambda})$. Then $X \cap ran (f_{\vec y}) \notin I_{\kappa, \pi}$.   

\medskip

{\bf Proof of Claim 1.} Given $t \in P_\kappa (\pi)$, set $b = \bigcup_{\beta \in t} y_\beta$. There must be $c \in P_\kappa (\lambda)$ such that $b \subseteq c$ and $f_{\vec y} (c) \in X$. Clearly $t \subseteq f_{\vec y} (c)$, which completes the proof of the claim. 

\medskip

{\bf Claim 2.}  Suppose that $X \cap ran (f_{\vec y}) \notin I_{\kappa, \pi}$. Then $X \notin f_{\vec y} (I_{\kappa, \lambda})$. 

\medskip

{\bf Proof of Claim 2.} Given $b \in P_\kappa (\lambda)$, we may find $c \in P_\kappa (\lambda)$ such that $f_{\vec y} (c) \in X$ and $b \subseteq f_{\vec y} (c)$. Then by Fact 3.11 (a), $b \subseteq f_{\vec y} (c) \cap \lambda = c$, which completes the proof of the claim and that of (ii). 
\hfill$\square$

\begin{fact} \begin{enumerate}[\rm (i)]
\item {\rm(\cite{Weaksat})} $\mathcal{A}_{\kappa,\lambda} (2,\lambda)$ holds, and in fact there is a $(2, \lambda,\pi)$-sequence $\vec y = \langle y_\beta : \beta < \pi\rangle$  with the property that $\vert \{y_\beta \cap a : \beta < \pi\} \vert < \kappa$ for all $a \in P_\kappa (\lambda)$.
\item {\rm(\cite{Weaksat})} If $\mathcal{A}_{\kappa,\lambda} (\tau,\pi)$ holds, then $\pi \leq \cov (\lambda, \kappa, \tau, 2)$.
\end{enumerate}  
\end{fact}

\medskip

Under certain cardinality assumptions, the existence of a pseudo-Kurepa family entails that of a Kurepa family of the same size.

\medskip

\begin{fact} {\rm(\cite{MPS2})} Suppose that $\mathcal{A}_{\kappa,\lambda} (\sigma^+,\lambda)$ holds, where $\lambda < \pi$ and $\sigma$ is a cardinal less than $\kappa$. Then there is $W \subseteq [\lambda^{< \sigma}]^{\cf (\sigma)}$ with $\vert W \vert = \pi$  such that $\vert \{w \in W : \vert w \cap x \vert = \cf (\sigma)\} \vert < \kappa$ for all $x \in P_\kappa (\lambda^{< \sigma})$.
\end{fact}

\begin{fact} {\rm(\cite{Todor})} Suppose that $\mathcal{A}_{\kappa,\lambda} (\sigma^+,\lambda)$ holds, where $\lambda < \pi$ and $\sigma$ is a cardinal less than $\kappa$. Suppose futher that $\lambda^{< \sigma} = \lambda$, and $\rho^{< \cf (\sigma)} < \kappa$ for any cardinal $\rho < \kappa$. Then there is a $((\cf (\sigma))^+, \lambda,\pi)$-sequence $\vec y = \langle y_\beta : \beta < \pi \rangle$  with the property that $\vert \{y_\beta \cap a : \beta < \pi\} \vert < \kappa$ for all $a \in P_\kappa (\lambda)$.
\end{fact}

{\bf Proof.} Use Fact 3.16.
\hfill$\square$  

\medskip

The existence of a pseudo-Kurepa family can be reformulated in several ways (for more on this, see \cite{MPS2} and \cite{Diamondstar}).

\medskip

\begin{fact} {\rm(\cite{Todor})} The following are equivalent :
\begin{enumerate}[\rm (i)]
\item ${\mathcal A}_{\kappa,\lambda} (\kappa, \pi)$.
\item There is $g : P_\kappa (\lambda) \rightarrow P_\kappa (\pi)$ such that (a) $g (a) \subseteq g (b)$ whenever $a \subseteq b$, and (b) $ran (g) \in I^+_{\kappa, \pi}$. 
\end{enumerate}
\end{fact}

{\bf Proof.} 

(i) $\rightarrow$ (ii) : By Observation 3.14.

\medskip

(ii) $\rightarrow$ (i) : By Fact 3.15 (i), $\mathcal{A}_{\kappa,\lambda} (\kappa,\lambda)$ holds, so we can assume that $\pi > \lambda$. Now given $g$ as in (ii), select a $(\kappa, \lambda,\pi)$-sequence $\vec y = \langle y_\beta : \beta < \pi \rangle$ so that for any $\beta \in \pi \setminus \lambda$, $\beta \in g (y_\beta)$.

\medskip

{\bf Claim 1.}  Let $a \in P_\kappa (\lambda)$. Then $\vert \{\beta \in \pi : y_\beta \subseteq a\} \vert < \kappa$.   

\medskip

{\bf Proof of Claim 1.}  It suffices to show that $\vert \{\beta \in \pi \setminus \lambda : y_\beta \subseteq a\} \vert < \kappa$, which follows from the fact that $\vert g (a) \vert < \kappa$, since for any $\beta \in \pi \setminus \lambda$ such that $y_\beta \subseteq a$, we have that $g (y_\beta) \subseteq g (a)$, and therefore $\beta \in g (a)$. Thus the proof of the claim is complete. 

\medskip

{\bf Claim 2.}  $\vert \{y_\beta : \beta \in  \pi\} \vert = \pi$.   

\medskip

{\bf Proof of Claim 2.} By Claim 1, the function $\beta \mapsto y_\beta$ is $< \kappa$-to-one, so

\centerline{$\pi \leq \max (\kappa, \vert \{y_\beta : \beta \in  \pi\} \vert)$.}

Since $\pi > \kappa$, the desired assertion follows, which completes the proof of the claim. 

\medskip

That ${\mathcal A}_{\kappa,\lambda} (\kappa, \pi)$ holds easily follows from Claims 1 and 2.
\hfill$\square$

\medskip

The next reformulation is in terms of a combinatorial principle of Rinot \cite{Rinot}.

\medskip

\begin{Obs}  Suppose  that $\pi$ is regular and greater than $\lambda$. Then for any cardinal $\sigma < \kappa$, the following are equivalent :
\begin{enumerate}[\rm (i)]
\item ${\mathcal A}_{\kappa,\lambda} (\sigma^+, \pi)$ holds.
\item There exists $d : \sigma \times \pi \rightarrow \lambda$ with the following property : for any $\alpha \in E^\pi_\kappa$ and any cofinal subset $B$ of $\alpha$, there are $\eta < \sigma$ and a cofinal subset $B'$ of $B$ such that the function $\beta \mapsto d (\eta, \beta)$ is one-to-one on $B'$.
\item There exists $d : \sigma \times \pi \rightarrow \lambda$ with the following property : for any $\alpha \in E^\pi_\kappa$ and any cofinal subset $B$ of $\alpha$, there are $\eta < \sigma$ and a subset $B'$ of $B$ with $\sup B' = \alpha$ such that the function $\beta \mapsto d (\eta, \beta)$ is $< \kappa$-to-one on $B'$.
\end{enumerate}
\end{Obs}

{\bf Proof.}

\hskip0,2cm  (i) $\rightarrow$ (ii) : Given an ${\mathcal A}_{\kappa,\lambda} (\sigma^+,\pi)$-sequence $\langle y_\beta : \beta < \pi \rangle$,  pick a bijection $f_\beta : \vert y_\beta \vert \rightarrow y_\beta$ for each $\beta \in \pi$, and define $d : \sigma \times \pi \rightarrow \lambda$ so that for any $\xi \in \pi$ and any $\eta < \vert y_\beta \vert$, $d (\eta, \xi)$ equals $f_\xi (\eta)$. Now fix $\alpha \in E^\pi_\kappa$. Select a cofinal subset $B$ of $\alpha$ of order-type $\kappa$. Notice that $\vert \bigcup_{\beta \in B} y_\beta \vert = \kappa$. Inductively define $\beta_i \in B$ and $\gamma_i \in y_{\beta_i}$ for $i < \kappa$ so that 
\begin{itemize}
\item $\beta_j < \beta_i$ whenever $j < i < \kappa$.
\item $\gamma_i \notin \bigcup_{j < i} y_{\beta_j}$ for all $i < \kappa$.
\end{itemize}
For $i < \kappa$, let $\gamma_i = f_{\beta_i} (\eta_i)$. There must be $\eta < \sigma$ and $e \in [\kappa]^\kappa$ such that $\eta_i = \eta$ for all $i \in e$. Then the function taking $i$ to $d (\eta, \beta_i)$ is one-to-one on $e$.

\hskip0,2cm  (ii) $\rightarrow$ (iii) : Trivial. 

\hskip0,2cm  (iii) $\rightarrow$ (i) : Let $d : \sigma \times \pi \rightarrow \lambda$ be as in (iii). For $\beta \in \pi \setminus \lambda$, put $y_\beta = \{d (\eta, \beta) : \eta < \sigma\}$. Now given $e \in [\pi \setminus \lambda]^\kappa$, select $B \subseteq e$ with $o.t. (B) = \kappa$, and put $\alpha = \sup e$. We may find $\eta \in \sigma$, and $B' \subseteq B$ with $\sup B' = \alpha$ such that the function taking $\beta$ to $d (\eta, \beta)$ is $< \kappa$-to-one on $B'$. Then clearly $\vert \bigcup_{\beta \in B'} y_\beta \vert = \kappa$, and consequently $\vert \bigcup_{\beta \in e} y_\beta \vert = \kappa$.
\hfill$\square$

\bigskip

\subsection{$A_\kappa (\vec y)$}

\bigskip

Let us show that if $\vec y = \langle y_\delta : \delta < \pi \rangle$ is an ${\mathcal A}_{\kappa,\lambda} (\kappa,\pi)$-sequence, then any $\kappa$-complete ideal $K$ on $P_\kappa (\pi)$ extending $I_{\kappa, \pi} \vert (ran (f_{\vec y}))$ is isomorphic to its projection.

\medskip

\begin{Def}  Given a $\kappa$-complete ideal $J$ on a set $X$, we let $\overline{cof} (J)$ denote the least size of any $Q \subseteq J$ with the property that for any $B \in J$, there is $q \in P_\kappa (Q)$ with $B \subseteq \bigcup q$.
\end{Def}

\begin{fact} {\rm (\cite{Secret})} Let $H$ be a fine ideal on $P_\kappa (\lambda)$, and $\vec y = \langle y_\delta : \delta < \pi \rangle$ be a $(\kappa, \lambda, \pi)$-sequence such that $A_\kappa (\vec y) \cap \Delta_\kappa (\vec y) \in (f_{\vec y} (H))^\ast$. Further let $K$ be a $\kappa$-complete ideal on $P_\kappa (\pi)$ extending $f_{\vec y} (H)$. Then $K = f_{\vec y} (p_\pi (K))$.
\end{fact} 

\begin{Obs} Let  $\vec y = \langle y_\beta : \beta < \pi \rangle$ be an ${\mathcal A}_{\kappa,\lambda} (\kappa,\pi)$-sequence, and $K$ be a $\kappa$-complete, fine ideal on $P_\kappa (\pi)$ with $ran (f_{\vec y}) \in K^\ast$. Then the following hold :
\begin{enumerate}[\rm (i)]
\item $K = f_{\vec y} (p_\pi (K))$.
\item Assuming that $K$ is normal, the following hold :
\begin{enumerate}[\rm (a)]
\item Let $C_\beta \in (p_\pi (K))^\ast$ for $\beta < \pi$. Then 

\centerline{$\{a \in P_\kappa (\lambda): \forall \beta \in f_{\vec y} (a) (a \in C_\beta)\} \in (p_\pi (K))^\ast$.}

\item Suppose that $J \subseteq p_\pi (K))$ is a $\kappa$-complete, fine ideal on $P_\kappa (\lambda)$ with $\overline{cof} (J) \leq \pi$. Then $J \vert C = I_{\kappa, \lambda} \vert C$ for some $C \in (p_\pi (K))^\ast$.
\end{enumerate}
\end{enumerate}
\end{Obs}

{\bf Proof.} (i) : Set $H = I_{\kappa, \lambda}$. Then by Fact 3.11 and Observation 3.14, $f_{\vec y} (H) = I_{\kappa, \pi} \vert (A_\kappa (\vec y) \cap \Delta_\kappa (\vec y)) = I_{\kappa, \pi} \vert ran (f_{\vec y})$. Now apply Fact 3.21.

\medskip

(ii) : By Lemma 2.1 and Proposition 2.2 in \cite{Norm}.

\hfill$\square$

\begin{Cor} Suppose that there exists an ${\mathcal A}_{\kappa,\lambda} (\kappa,\pi)$-sequence $\vec y$ such that $A_\kappa (\vec y) \in NS_{\kappa, \pi}^+$, where $\pi = 2^\lambda$. Then $NS_{\kappa, \lambda} \vert C = I_{\kappa, \lambda} \vert C$ for some $C \in NS_{\kappa, \lambda}^+$.
\end{Cor}

{\bf Proof.} Apply Observation 3.22 with $K = NS_{\kappa, \pi} \vert A_\kappa (\vec y)$.
\hfill$\square$

\medskip

Note that if $NS_{\kappa, \lambda} \vert C = I_{\kappa, \lambda} \vert C$ for some $C$, then as shown in \cite{MPS1}, $cof (NS_{\kappa, \lambda}) = u (\kappa, \lambda)$, and therefore by a result of \cite{MS}, $NS_{\kappa, \lambda}$ is nowhere precipitous.

\medskip

\begin{Obs} Let  $\vec y = \langle y_\beta : \beta < \pi \rangle$ be an ${\mathcal A}_{\kappa,\lambda} (\kappa,\pi)$-sequence, and $J$ be a fine ideal on $P_\kappa (\lambda)$. Then $p_\pi (f_{\vec y} (J)) = J$.
\end{Obs}

{\bf Proof.} Given $B \subseteq P_\kappa (\lambda)$,

$B \in p_\pi (f_{\vec y} (J))$ 

iff $\{x \in P_\kappa (\pi) : x \cap \lambda \in B\} \in f_{\vec y} (J)$

iff $\{b \in P_\kappa (\lambda) : f_{\vec y} (b) \cap \lambda \in B\} \in J$ 

iff $B \in J$.
\hfill$\square$

\begin{fact} {\rm (\cite{Secret})} \begin{enumerate}[\rm (i)]
\item Let $K$ be a fine ideal on $P_\kappa (\pi)$, and $\vec y$ be a $(\kappa, \lambda, \pi)$-sequence such that $A_\kappa (\vec y) \cap \Delta_\kappa (\vec y) \in K^\ast$. Then for any $A \in (p_\pi (K))^+$, there is $W \in K^+$ such that $p_\pi (K) \vert A = p_\pi (K \vert W)$. 
\item Let $J$ be a fine ideal on $P_\kappa (\lambda)$, and $\vec y$ be a $(\kappa, \lambda, \pi)$-sequence such that $A_\kappa (\vec y) \cap \Delta_\kappa (\vec y) \in (f_{\vec y} (J))^\ast$. Then for any $W \in (f_{\vec y} (J))^+$, there is $A \in J^+$ such that $f_{\vec y} (J) \vert W = f_{\vec y} (J \vert A)$. 
\end{enumerate}
\end{fact} 

\begin{Obs} Let $\vec y = \langle y_\delta : \delta < \pi \rangle$ be an ${\mathcal A}_{\kappa,\lambda} (\kappa,\pi)$-sequence. Then the following hold : 
\begin{enumerate}[\rm (i)]
\item  For any $A \in I_{\kappa, \lambda}^+$, there is $W \in (I_{\kappa, \pi} \vert ran (f_{\vec y}))^+$ such that (a) $I_{\kappa, \lambda} \vert A = p_\pi ((I_{\kappa, \pi} \vert ran ( f_{\vec y})) \vert W)$, and (b) $f_{\vec y} (I_{\kappa, \lambda} \vert A) = (I_{\kappa, \pi} \vert ran ( f_{\vec y})) \vert W$. 
\item For any $W \in (I_{\kappa, \pi} \vert ran ( f_{\vec y}))^+$, there is $A \in I_{\kappa, \lambda}^+$ such that (a) $(I_{\kappa, \pi} \vert ran ( f_{\vec y})) \vert W = f_{\vec y} (I_{\kappa, \lambda} \vert A)$, and (b) $p_\pi ((I_{\kappa, \pi} \vert ran ( f_{\vec y})) \vert W) = I_{\kappa, \lambda} \vert A$.
\end{enumerate}
\end{Obs}

{\bf Proof.} By Fact 3.25 and Observations 3.22 and 3.24.
\hfill$\square$

\medskip

Suppose for instance that $\lambda$ is a strong limit cardinal of cofinality less than $\kappa$, and $\vec y$ is an ${\mathcal A}_{\kappa,\lambda} (\kappa,\pi)$-sequence, where $\pi = 2^\lambda$. Then by a result of Shelah \cite{She02}, $NS_{\kappa, \lambda}$ is a restriction of $I_{\kappa, \lambda}$, so by Observation 3.26,  for any $S \in NS_{\kappa, \lambda}^+$, there is $T \in I_{\kappa, \pi}^+$ such that $f_{\vec y} (NS_{\kappa, \lambda} \vert S) = I_{\kappa, \pi} \vert T$. 

\bigskip

\subsection{Existence of shuttles : the case $\cf (\lambda) < \mu$}

\bigskip

We will now investigate conditions for the existence of shuttles. We start with the easier case, that is when $\cf (\lambda) < \mu$.

\medskip

\begin{Pro} Let  $\vec y = \langle y_\beta : \beta < \pi \rangle$ be an ${\mathcal A}_{\kappa,\lambda} (\mu,\pi)$-sequence. Then $f_{\vec y}$ is a $(\mu, \kappa, \lambda, \pi)$-shuttle.
\end{Pro}

{\bf Proof.} Put $\chi = f_{\vec y}$. 
Clauses (b) and (d) of Definition 3.1 are easily seen to be satisfied. For clause (a) (respectively, (c)), use Fact 3.11 (i) (respectively, Observation 3.14 (i)). 
\hfill$\square$

\medskip

There are situations in which ${\mathcal A}_{\kappa,\lambda} (\mu,\pi)$ follows from ${\mathcal A}_{\kappa,\lambda} (\kappa,\pi)$.
\medskip

\begin{fact}\begin{enumerate}[\rm (i)] 
\item  {\rm(\cite{Norm})}  Suppose that 
\begin{itemize}
\item ${\mathcal A}_{\kappa,\lambda} (\kappa,\pi)$ hold ;
\item $\cf (\lambda) < \mu$ ;
\item $\cov (\sigma, \kappa, \kappa, \mu) \leq \lambda$ for any cardinal $\sigma$ with $\kappa \leq \sigma < \lambda$. 
\end{itemize}
Then ${\mathcal A}_{\kappa,\lambda} (\mu,\pi)$ holds. 
\item {\rm (\cite{Secret})} Suppose that 
\begin{itemize}
\item either $\lambda < \kappa^{+ \mu}$,
\item or $\cf (\lambda) < \mu$ and SSH holds.
\end{itemize}
Then ${\mathcal A}_{\kappa,\lambda} (\kappa,\pi)$ and ${\mathcal A}_{\kappa,\lambda} (\mu,\pi)$ are equivalent. 
\end{enumerate}
\end{fact}

\begin{Def} Suppose that $\tcf (\prod A/I) = \pi$, where $A$ is an infinite set of regular infinite cardinals such that $\vert A \vert  < \min (A)$, and $I$ an ideal on $A$, and that $\vec{f} = \langle f_\xi : \xi < \pi \rangle$ is an $<_I$-increasing, cofinal sequence in $(\prod A, <_I)$. An infinite limit ordinal $\delta < \pi$  is a {\it good point} for $\vec f$ if we may find a cofinal subset $X$ of $\delta$,  and $Z_\xi\in I$ for $\xi\in X$ such that $f_\beta (i) < f_\xi (i)$  whenever $\beta < \xi$  are in $X$ and $i \in \nu \setminus (Z_\beta \cup Z_\xi)$. 
\end{Def}

\medskip

The following is well-known (see e.g. \cite[Lemma 2.5]{AM}).

\begin{fact} Let $\delta < \pi$  be a good point for $\vec f$, and $e$ be a cofinal subset of $\delta$ of order-type $\cf (\delta)$. Then  there are $X \in [e]^{\cf (\delta)}$,  and $Z_\xi \in I$ for $\xi \in X$ such that $f_\beta (a) < f_\xi (a)$  whenever $\beta < \xi$  are in $X$ and $a \in A \setminus (Z_\beta \cup Z_\xi)$. 
\end{fact}

\begin{fact}  {\rm(\cite{Good}, \cite{Weaksat})} Let $\theta, \pi, A$ and $I$ be such that
\begin{itemize}
\item $\theta$ and $\pi$ are two cardinals such that $\kappa < \theta < \pi$ ;
\item $A$ is a set of regular cardinals smaller than $\theta$ such that $\vert A \vert < \kappa$ and $\sup A = \theta$ ;
\item $I$ is an ideal on $A$ with $\{A \cap a : a \in A \} \subseteq I$ ;
\item  $\pi = \tcf (\prod A /I )$.
\end{itemize}
Further let  $\vec{f} = \langle f_\alpha : \alpha < \pi \rangle$ be an increasing, cofinal sequence in $(\prod A, <_I)$. Then the following hold :
\begin{enumerate}[\rm (i)]
\item There is a closed unbounded subset $C_{\vec{f}}$ of $\pi$, consisting of infinite limit ordinals, with the property that any $\delta$ in $C_{\vec{f}}$ satisfying one of the following conditions, where $\rho$ denotes the largest limit cardinal less than or equal to $\cf (\delta)$, is a good point for $\vec f$ :
\begin{enumerate}[\rm (a)]
\item $(\max (\rho, \vert A \vert))^{+3} < \cf (\delta)$.
\item $\rho^{\vert A \vert} < \cf (\delta)$.
\item $\vert A \vert < \cf (\rho)$.
\item $\vert A \vert < \rho$  and $I$ is $(\cf (\rho))^+$-complete.
\item $\cf (\rho) \leq \vert A \vert < \rho$ and $\pp (\rho) < \cf (\delta)$. 
\end{enumerate}
\item Suppose that $\kappa$ is a successor, $\theta \leq \lambda < \pi$, and there is a closed unbounded subset $C$ of $\pi$ such that every $\delta \in C$ of cofinality $\kappa$ is a good point for ${\vec f}$. Then ${\mathcal A}_{\kappa,\lambda} (\vert A \vert ^+, \pi)$ holds.
\item Suppose that $\kappa$  is weakly inaccessible and $\theta \leq \lambda < \pi$. Then ${\mathcal A}_{\kappa,\lambda} (\vert A \vert ^+, \pi)$ holds.
\end{enumerate}
\end{fact}

\begin{Pro} Assume that $\kappa$  is weakly inaccessible and $\cf (\lambda) < \mu$. Then  ${\mathcal A}_{\kappa,\lambda} (\mu, \lambda^+)$ holds.
\end{Pro}

{\bf Proof.} By Facts 1.4 and 3.31.
\hfill$\square$

\begin{Pro} Suppose that SSH holds, $\cf (\lambda) < \mu$, $\kappa$ is the successor of a cardinal $\nu$, and $N\mu$-$S_{\kappa,\lambda}$ and $N\mu$-$S_{\kappa,u (\kappa, \lambda)}$ are not isomorphic. Then (by Propositions 3.7 and 3.27, ${\mathcal A}_{\kappa,\lambda} ((\cf (\lambda))^+, u (\kappa,\lambda))$ does not hold, and) one of the following holds :
\begin{enumerate}[\rm (a)] 
\item $\mu = \nu = (\cf (\lambda))^+$ and $2^{\cf (\lambda)} \geq \kappa$.
\item $\mu = \nu = (\cf (\lambda))^{++}$ and $2^{\cf (\lambda)} \geq \kappa$.
\item $\mu = (\cf (\lambda))^+$ and $\nu = (\cf (\lambda))^{++}$ and $2^{\cf (\lambda)} \geq \kappa$.
\item $\cf (\nu) = \cf (\lambda) < \nu$.
\end{enumerate}
\end{Pro}

{\bf Proof.} Let $\rho$ denote the largest limit cardinal less than $\kappa$. Then by Facts 1.4 and 3.31, the following hold :
\begin{itemize}
\item $\cf (\lambda) \geq \cf (\rho)$.
\item If $\cf (\lambda) < \rho$, then $\rho$ is singular, and moreover (1) $\rho^+ \geq \kappa$, so $\rho = \nu$, and (2) $\cf (\lambda) \leq \cf (\rho)$, so $\cf (\lambda) = \cf (\rho) = \cf (\nu)$.
\item If $\cf (\lambda) \geq \rho$, then $2^{\cf (\lambda)} \geq \kappa$, and moreover $\kappa \leq (\cf (\lambda))^{+ 3}$, so $\cf (\lambda) \leq \nu \leq (\cf (\lambda))^{+ +}$.
\end{itemize}
\hfill$\square$

\begin{Cor} Suppose that SSH holds, $\cf (\lambda) < \mu$, $2^{\cf (\lambda)} < \kappa$, $\kappa$ is the successor of a cardinal $\nu$, and $N\mu$-$S_{\kappa,\lambda}$ and $N\mu$-$S_{\kappa,u (\kappa, \lambda)}$ are not isomorphic. Then $\cf (\nu) = \cf (\lambda) < \nu$.
\end{Cor}


\bigskip

\subsection{More on the case $\cf (\lambda) < \mu$}

\bigskip

We now turn to $(\mu, \kappa)$-normal ideals. Let us start by recalling the following.

\medskip

\begin{fact} {\rm (\cite{Secret})} Let $\vec y = \langle y_\beta : \beta < \pi \rangle$ be an ${\mathcal A}_{\kappa,\lambda} (\mu,\pi)$-sequence. Then $A_\kappa (\vec y)  \in NS^\ast_{\mu, \kappa, \pi}$.
\end{fact}

\medskip

Next we show that if ${\mathcal A}_{\kappa,\lambda} (\mu,\pi)$ holds, then any $(\mu, \kappa)$-normal ideal $J$ on $P_\kappa (\lambda)$ is isomorphic to a $(\mu, \kappa)$-normal ideal $K$ on $P_\kappa (\pi)$ (of which by Observation 3.24 it is the projection).

\medskip

\begin{Obs} Let  $\vec y = \langle y_\beta : \beta < \pi \rangle$ be an ${\mathcal A}_{\kappa,\lambda} (\mu,\pi)$-sequence, and $J$ be a $(\mu, \kappa)$-normal ideal on $P_\kappa (\lambda)$. Then the following hold :
\begin{enumerate}[\rm (i)] 
\item $f_{\vec y} (J)$ is a $(\mu, \kappa)$-normal ideal on $P_\kappa (\pi)$ (which by Fact 2.7, Observation 2.9 and Observation 3.14 is isomorphic to $J$).
\item For any $A \in J^+$, there is $W \in (f_{\vec y} (J))^+$ such that $f_{\vec y} (J \vert A) = f_{\vec y} (J) \vert W$.
\end{enumerate}
\end{Obs}

{\bf Proof.} (i) : By Fact 2.4 and Observation 3.14, $f_{\vec y} (J)$ is a fine, $\kappa$-complete ideal on $P_\kappa (\pi)$. Now let $X \in (f_{\vec y} (J))^+$, and  $g : X \rightarrow P_\mu (\pi)$ with the property that $g (x) \subseteq x$ for every $x \in X$. Define  $h: f^{- 1}_{\vec y} (X) \rightarrow P_\mu (\lambda)$ by $h (a) = \bigcup \{y_\alpha : \alpha \in g (f_{\vec y} (a)) \}$. We may find $w \in P_\kappa ( \lambda)$ and $A \in J^+ \cap P (f^{- 1}_{\vec y} (X))$ such that $h (a) \subseteq w$ for all $a \in A$. Set $T = f_{\vec y}`` A$ and $z = f_{\vec y} (w)$. Then $T \in (f_{\vec y} (J))^+ \cap P (X)$, and moreover $g (x) \subseteq z$ for every $x \in T$.

\medskip

(ii) : Fix $A \in J^+$. Set $K = f_{\vec y} (J)$. Then by Observation 3.24, $p_\pi (K) = J$. By (i) and Fact 3.35, $A_\kappa (\vec y) \in K^\ast$, so by Fact 3.25 (i), there is $W \in K^+$ such that $J \vert A = p_\pi (K \vert W)$. Finally, by Observation 3.22, $f_{\vec y} (J \vert A) = K \vert W$.
\hfill$\square$


\medskip

It is also true that if ${\mathcal A}_{\kappa,\lambda} (\mu,\pi)$ holds, then each $(\mu, \kappa)$-normal ideal $K$ on $P_\kappa (\pi)$ is isomorphic to a $(\mu, \kappa)$-normal ideal $J$ on $P_\kappa (\lambda)$. 

\medskip

\begin{Obs} Let  $\vec y = \langle y_\beta : \beta < \pi \rangle$ be an ${\mathcal A}_{\kappa,\lambda} (\mu,\pi)$-sequence, and $K$ be a $(\mu, \kappa)$-normal ideal on $P_\kappa (\pi)$. Then the following hold : 
 \begin{enumerate}[\rm (i)] 
\item $K = f_{\vec y} (p_\pi (K))$ (where $p_\pi (K)$ is $(\mu, \kappa)$-normal by Fact 2.22 (iv)).
\item For any $W \in K^+$, there is $A \in (p_\pi (K))^+$ such that  $K \vert W = f_{\vec y} (p_\pi (K) \vert A)$.
\end{enumerate}
\end{Obs}

{\bf Proof.} By Observation 3.22 and Facts 3.25 (ii) and 3.35.
\hfill$\square$

\medskip

\begin{Rmk} Suppose that SSH holds, $\cf (\lambda) < \mu$ and $\pi = \lambda^+ (= u (\kappa, \lambda))$. Then by results of \cite{Secret}, $NS_{\mu, \kappa, \pi} = SNS_{\kappa, \pi} \vert S$ for some $S$, and moreover any normal extension of $NS_{\mu, \kappa, \pi}$ is $(\mu, \kappa)$-normal.
\end{Rmk}

\medskip

Let us now consider specific $(\mu, \kappa)$-normal ideals.

\medskip

\begin{Pro} Let  $\vec y = \langle y_\beta : \beta < \pi \rangle$ be an ${\mathcal A}_{\kappa,\lambda} (\mu,\pi)$-sequence. Then $f_{\vec y} (NS_{\mu, \kappa, \lambda}) = NS_{\mu, \kappa, \pi}$.
\end{Pro}

{\bf Proof.}
By Fact 2.22 (iv), $NS_{\mu, \kappa, \lambda} \subseteq p_\pi (NS_{\mu, \kappa, \pi})$, so by Observation 3.37, $f_{\vec y} (NS_{\mu, \kappa, \lambda}) \subseteq f_{\vec y} (p_\pi (NS_{\mu, \kappa, \pi})) = NS_{\mu, \kappa, \pi}$. On the other hand by Observation 3.36, $NS_{\mu, \kappa, \pi} \subseteq f_{\vec y} (NS_{\mu, \kappa, \lambda})$.
\hfill$\square$

\medskip

Now the corresponding result for the game ideal.

\medskip

\begin{Pro} Let  $\vec y = \langle y_\beta : \beta < \pi \rangle$ be an ${\mathcal A}_{\kappa,\lambda} (\mu,\pi)$-sequence. Then $f_{\vec y} (NG^\mu_{\kappa, \lambda}) = NG^\mu_{\kappa, \pi}$.
\end{Pro}

{\bf Proof.} By Fact 2.25 and Observation 3.37.
\hfill$\square$

\bigskip

\subsection{Existence of shuttles : the case $\mu < \cf (\lambda) < \kappa$}

\bigskip

The case $\mu < \cf (\lambda) < \kappa$ will be dealt with by appealing to pseudo-Kurepa families of a special kind.

\medskip

The following freeness principle is closely related to several principles of Shelah (see e.g. \cite[Definition 1.10]{She1008}).

\medskip

\begin{Def} Let $\sigma < \kappa$ be an infinite cardinal, and $I$ be an ideal on $\sigma$. 
An ${\mathcal F}_{\kappa,\lambda}^I (\sigma^+, \pi)${\it - sequence}
is a $(\sigma^+,\lambda,\pi)$-sequence $\vec y = \langle y_\beta : \beta < \pi \rangle$ for which one may find a sequence $\vec f = \langle f_\beta : \lambda \leq \beta < \pi \rangle$ such that :
\begin{itemize}
\item For each $\beta \in \pi \setminus \lambda$, $f_\beta$ is a function from $\sigma$ to $ \lambda$.
\item  For any $e \in [\pi \setminus \lambda]^\kappa$, there are $b \in [e]^\kappa$ and $g : b \rightarrow I$ such that $f_\alpha (i) \not= f_\beta (i)$  whenever $\alpha < \beta$ are in $b$ and $i$ is in $\sigma \setminus (g(\alpha) \cup g(\beta))$.
\item For $\lambda \leq \beta < \pi$, $y_\beta = ran (f_\beta)$.
\end{itemize}
\end{Def} 




\begin{Obs} Let $\vec y = \langle y_\beta : \beta < \pi\rangle$ be an ${\mathcal F}_{\kappa,\lambda}^I   (\sigma^+, \pi)$- sequence, where $\sigma < \kappa$ is an infinite cardinal and $I$ an ideal on $\sigma$, as witnessed by $\vec f = \langle f_\beta : \lambda \leq \beta < \pi \rangle$. Then for any $v \in P_\kappa (\lambda)$,

\centerline{$\vert \{ \beta \in \pi \setminus \lambda : \{i < \sigma : f_\beta (i) \in v\} \in I^+\} \vert < \kappa$.} 
\end{Obs}

{\bf Proof.}  Fix $v \in P_\kappa (\lambda)$. Suppose toward a contradiction that there is $e \in [\pi \setminus \lambda]^\kappa$ such that $\{i < \sigma : f_\beta (i) \in v\} \in I^+$ for all $\beta \in e$. We may find $b \in [e]^\kappa$ and $g : b \rightarrow I$ such that $f_\alpha (i) \not= f_\beta (i)$  whenever $\alpha < \beta$ are in $b$ and $i$ is in $\sigma \setminus (g(\alpha) \cup g(\beta))$. For $\gamma \in b$, pick $i_\gamma \in \sigma \setminus g (\gamma)$ with $f_\gamma (i_\gamma) \in v$. There must be $i < \sigma$ and $c \in [b]^\kappa$ such that $i_\gamma = i$ for all $\gamma \in c$. But then  $f_\alpha (i) \not= f_\beta (i)$  whenever $\alpha < \beta$ are in $c$. This is a contradiction since $\vert v \vert < \kappa$.
\hfill$\square$

\begin{Cor}  Let $\vec y = \langle y_\beta : \beta < \pi\rangle$ be an ${\mathcal F}_{\kappa,\lambda}^I   (\sigma^+, \pi)$- sequence, where $\sigma < \kappa$ is an infinite cardinal, and $I$ an ideal on $\sigma$. Then $\vec y$ is an ${\mathcal A}_{\kappa,\lambda} (\sigma^+,\pi)$-sequence.
\end{Cor}

\begin{Obs} Suppose that there is an ${\mathcal F}_{\kappa,\lambda}^I   (\sigma^+, \pi)$- sequence $\vec y = \langle y_\beta : \beta < \pi\rangle$, where $\sigma < \kappa$ is an infinite cardinal and $I$ a $\mu^+$-complete ideal on $\sigma$. Then $\pi \leq \cov (\lambda, \kappa, \kappa, \mu^+)$.
\end{Obs}

{\bf Proof.}  By the proof of Corollary 11.4 in \cite{Secret}.
\hfill$\square$

\begin{Pro} Let $\vec y = \langle y_\beta : \beta < \pi\rangle$ be
an ${\mathcal F}_{\kappa,\lambda}^I   (\sigma^+, \pi)$- sequence, where $\sigma < \kappa$ is an infinite cardinal and $I$ a $\mu^+$-complete ideal on $\sigma$, as witnessed by $\vec f = \langle f_\beta : \lambda \leq \beta < \pi \rangle$. Define $\chi : P_\kappa (\lambda) \rightarrow P_\kappa (\pi)$ by 

\centerline{$\chi (a) = a \cup \{ \beta \in \pi \setminus \lambda : \{i < \sigma : f_\beta (i) \in a \} \in I^+\}$.}

Then the following hold :
\begin{enumerate}[\rm (i)]
\item For any $b \in P_\kappa (\lambda)$, $f_{\vec y} (b) \subseteq \chi (b)$.
\item For any $u : \mu \rightarrow P_\kappa (\lambda)$,  $\chi (\bigcup_{r < \mu} u (r)) \subseteq \bigcup_{r < \mu} \chi (u (r))$
\item $\chi$ is a $(\mu, \kappa, \lambda, \pi)$-shuttle.
\item $A_\kappa (\vec y) \in (NG_{\kappa, \pi}^\mu)^\ast$.
\item $f_{\vec y} (NG^\mu_{\kappa, \lambda}) = NG^\mu_{\kappa, \pi}$.
\end{enumerate}
\end{Pro}

{\bf Proof.}  (i) : Immediate.

\medskip

(ii) : By the proof of Proposition 11.5 in \cite{Secret}.

\medskip

(iii) : Clauses (a) and (b) of Definition 3.1 are obviously verified. For clause (d), use (ii). For clause (c), fix $w \in P_\kappa (\pi)$. Proceeding inductively, define $x_r$ for $r < \mu$ as follows :
\begin{itemize}
\item $x_0 = (w \cap \lambda) \cup (\bigcup_{\beta \in w \setminus \lambda} ran (f_\beta))$.
\item $x_r = \bigcup_{s < r} x_s$ in case $r$ is an infinite limit ordinal.
\item $x_{r + 1} = x_r \cup\chi (x_r \cap \lambda)$.
\end{itemize}
Put $x = \bigcup_{r < \mu} x_r$. Clearly, $w \subseteq \chi (x_0) \subseteq x_1 \subseteq x$. Furthermore, 

\centerline{$\chi (x \cap \lambda) \subseteq \bigcup_{r < \mu} \chi (x_r \cap \lambda) \subseteq \bigcup_{r < \mu} x_{r + 1} \subseteq x$.}

Let us finally establish that $x \subseteq \chi (x \cap \lambda)$. To show this, fix $\gamma \in x$, and let $s$ be least such that $\gamma \in x_s$. If $s$ is a successor ordinal, say $s = r + 1$, then $\gamma$ lies in $\chi (x_r \cap \lambda)$ and hence in $\chi (x \cap \lambda)$. Otherwise, $s = 0$. But then $\gamma$ lies in $x \cap \lambda$ and therefore in $\chi (x \cap \lambda)$.

\medskip

(iv) : We define a winning strategy $\tau$ for player II in $G_{\kappa, \pi}^\mu (A_\kappa (\vec y))$ by : $\tau (x_0, x_1, ... , x_j) = x_j \cup \chi (x_j \cap \lambda)$.

\medskip

(v) : Set $K = NG^\mu_{\kappa, \pi}$. Then by (iv), Fact 2.25 (ii), Observation 3.22 and Corollary 3.43, $K =  f_{\vec y} (p_\pi (K)) = f_{\vec y} (NG^\mu_{\kappa, \lambda})$.
\hfill$\square$
 
\begin{Def}  A {\it beautiful} $(\mu, \kappa, \lambda, \pi)${\it-bridge} is an ${\mathcal A}_{\kappa,\lambda} (\kappa,\pi)$-sequence $\vec y = \langle y_\beta : \beta < \pi \rangle$ for which there is $\chi : P_\kappa (\lambda) \rightarrow P_\kappa (\pi)$ with the following properties : 
\begin{itemize}
\item Let $u : \mu \rightarrow P_\kappa (\lambda)$ be such that $\bigcup_{i < \mu} u (i) = \bigcup_{n \in e} u (n)$ for every $e \in [\mu]^\mu$. Then $f_{\vec y} (\bigcup_{i < \mu} u (i)) \subseteq \bigcup_{i < \mu} \chi (u (i))$. 
\item $\chi (a) \subseteq \chi (b)$ whenever $a \subseteq b$.
\item $\chi (\bigcup_{i < \mu} a_i) \subseteq \bigcup_{i < \mu} \chi (a_i)$ for any increasing sequence $\langle a_i : i < \mu \rangle$ in $(P_\kappa (\lambda), \subset)$. 
\end{itemize}
\end{Def}

\begin{Cor} Let $\vec y = \langle y_\beta : \beta < \pi\rangle$ be an ${\mathcal F}_{\kappa,\lambda}^I   (\sigma^+, \pi)$- sequence, where $\sigma < \kappa$ is an infinite cardinal and $I$ a $\mu^+$-complete ideal on $\sigma$, as witnessed by $\vec f = \langle f_\beta : \lambda \leq \beta < \pi \rangle$. Then $\vec y$ is a beautiful $(\mu, \kappa, \lambda, \pi)$-bridge.
\end{Cor}

\medskip

\begin{Rmk} The corresponding result for the case when $\cf (\lambda) < \mu$ also holds. Namely, as shown in \cite{Secret}, if $\vec y = \langle y_\beta : \beta < \pi \rangle$ is an ${\mathcal A}_{\kappa,\lambda} (\mu,\pi)$-sequence, then $\vec y$ is a beautiful $(\mu, \kappa, \lambda, \pi)$-bridge.
\end{Rmk}

\medskip

Beautiful bridges were used in \cite{Secret} to give an alternative characterization of  $N\mu$-$S_{\kappa,\lambda}$.

\begin{Def}  For a cardinal $\theta \geq \lambda$ and $h : \theta \rightarrow P_\kappa(\lambda)$, we let $X^h_{\mu, \kappa, \lambda, \theta}$ denote the set of all $x \in P_\kappa(\theta)$ for which there is $z \in [x]^\mu$ such that $x \cap \lambda \subseteq \bigcup_{\alpha \in t} h(\alpha)$ for all $t \in [z]^\mu$.
\end{Def}

\begin{fact} Assume that there exists a beautiful $(\mu, \kappa, \lambda, \pi)$-bridge $\vec y$, and $X^k_{\mu, \kappa, \lambda, \pi} \in SNS^+_{\kappa, \pi}$ for some $k : \pi \rightarrow P_\kappa(\lambda)$. Then there is $h : \pi \rightarrow P_\kappa(\lambda)$ such that $N\mu$-$S_{\kappa,\lambda} = p_\pi (SNS_{\kappa, \pi} \vert X^h_{\mu, \kappa, \lambda, \pi})$.
\end{fact}

\medskip

Let us finally give a condition for the existence of ${\mathcal F}_{\kappa,\lambda}^I (\sigma^+, \pi)$- sequences. 

\medskip

\begin{Obs} Suppose that $\tcf (\prod A/I) = \pi$, where $A$ is an infinite set of regular infinite cardinals with $\vert A \vert = \sigma < \min (A)$ and $\kappa < \sup A \leq \lambda < \pi$, and $I$ an ideal on $A$, and that $\vec{g} = \langle g_\xi : \xi < \pi \rangle$ is an $<_I$-increasing, cofinal sequence in $(\prod A, <_I)$. Suppose further that for some closed unbounded subset $C$ of $\pi$, every $\delta$ in $C \cap E^\pi_\kappa$  is a good point for $\vec g$. Then there exists an ${\mathcal F}_{\kappa,\lambda}^I (\sigma^+, \pi)$- sequence.
\end{Obs}

{\bf Proof.} Let $\langle \alpha_\beta : \lambda \leq \beta < \pi \rangle$ enumerate the elements of $C$ in increasing order. For $\lambda \leq \beta < \pi$, put $f_\beta = g_{\alpha_\beta}$. Given $e \in [\pi \setminus \lambda]^\kappa$, select a subset $e'$ of $e$ of order-type $\kappa$, and set $\delta = \sup \{ \alpha_\beta : \beta \in e' \}$. Then $\delta$ lies in $C \cap E^\pi_\kappa$, so it is a good point for $\vec g$. By Fact 3.30, we may find $X$ in $[e']^{\kappa}$ (and hence in $[e]^{\kappa}$), and $Z_\xi \in I$ for $\xi \in X$ such that $f_\beta (a) = g_{\alpha_\beta} (a) < g_{\alpha_\xi} (a) = f_\xi (a)$ whenever $\beta < \xi$  are in $X$ and $a \in A \setminus (Z_\beta \cup Z_\xi)$.
\hfill$\square$  


 

\bigskip

\subsection{Summary}
 
\bigskip

If $u (\kappa, \lambda)$ is regular and $\cf (\lambda) \not= \mu$, then under SSH, either $\lambda$ is regular and $u (\kappa, \lambda) = \lambda$, or $\cf (\lambda) \in \kappa \setminus \{ \mu\}$ and $u (\kappa, \lambda) = \lambda^+$. So Theorem 1.1 will follow from the following. 

\medskip

\begin{Pro} Suppose that one of the following holds :
\begin{enumerate}[\rm (i)]
\item $\cf (\lambda) < \mu$ and ${\mathcal A}_{\kappa,\lambda} (\mu, \lambda^+)$.
\item $\mu < \cf (\lambda) < \kappa$ and for some infinite set $A$ of regular infinite cardinals with $\vert A \vert < \min (\kappa, \min (A))$ and $\sup A = \lambda$, some $\mu^+$-complete ideal $I$ on $A$, and some $<_I$-increasing, cofinal sequence $\vec{g} = \langle g_\xi : \xi < \lambda^+ \rangle$ in $(\prod A, <_I)$, there is a closed unbounded subset $C$ of $\lambda^+$ such that every $\delta$ in $C \cap E^{\lambda^+}_\kappa$  is a good point for $\vec g$. 
\end{enumerate}
Then there exists a $(\mu, \kappa, \lambda, \lambda^+)$-shuttle (and hence by Proposition 3.7, $N\mu$-$S_{\kappa,\lambda}$ and $N\mu$-$S_{\kappa,\lambda^+}$ are isomorphic). 
\end{Pro}

{\bf Proof.} Case $\cf (\lambda) < \mu$ : By Proposition 3.27.

Case $\mu < \cf (\lambda) < \kappa$ :  By Proposition 3.45 and Observation 3.51.
\hfill$\square$

\bigskip

\section{Saturation}

\bigskip

In this section we look for cardinals $\rho$ such that $N\mu$-$S_{\kappa,\lambda}$ (and hence $NS_{\kappa, \lambda} \vert E^{\kappa, \lambda}_\mu$) is not $\rho$-saturated. Earlier papers on non-saturation properties of $NS_{\kappa, \lambda}$ include \cite{BurMats}, \cite{ForMag2} and \cite{Usuba}. In the case when $\lambda$ is regular, our approach is the usual one based on the function $f$ defined on $P_\kappa (\lambda)$ by $f (b) = \sup b$. Thus in this case our results on $P_\kappa (\lambda)$ 
will rely on known results on $\lambda$, but also take them as models. In the case when $\cf (\lambda) \in \kappa \setminus \{ \mu \}$, we use shuttles to deduce non-saturation properties on $P_\kappa (\lambda)$ from corresponding non-saturation properties on $P_\kappa (\lambda^+)$.

\bigskip

\subsection{$N\mu$-$S_{\kappa,\lambda}$ and the Sup function}

\bigskip

\begin{Def} We define $Sup : P_\kappa (\lambda) \rightarrow \lambda + 1$ by $Sup (a) = \sup a$.
\end{Def}

\begin{Obs} Suppose that $\lambda$ is regular. Then the following hold :
\begin{enumerate}[\rm (i)]
\item For any $S, T \subseteq \lambda$, $Sup^{- 1} (S \setminus T) = Sup^{- 1} (S) \setminus Sup^{- 1} (T)$.
\item $I_\lambda \subseteq Sup (I_{\kappa, \lambda})$.
\item $NS_\lambda \subseteq Sup (SNS_{\kappa, \lambda})$.
\end{enumerate}
\end{Obs}

{\bf Proof.} (i) and (ii) are immediate. For (iii), use Fact 2.20.
\hfill$\square$ 

\medskip

In the case when $\lambda$ is regular, our study of non-saturation properties of a given ideal $K$ on $P_\kappa (\lambda)$ will very much depend on properties of the associated ideal $Sup (K)$. The remainder of this subsection is devoted to $Sup (N\mu$-$S_{\kappa,\lambda})$ and $Sup (NG^\mu_{\kappa,\lambda})$.

\medskip

\begin{fact}  {\rm(\cite{Resemble}, \cite{Secret})} Suppose that $\lambda = \cf (\lambda)$, and $u (\kappa, \sigma) < \lambda$ for every cardinal $\sigma$ with $\kappa \leq \sigma < \lambda$. Then $Sup (N\mu$-$S_{\kappa,\lambda}) \subseteq NS_\lambda \vert E^\lambda_\mu$.     
\end{fact} 

\medskip

We will show that the conclusion of Fact 4.3 remains valid if $\lambda =\nu^+$ with $\cf (\nu) = \mu$ provided that $u (\kappa, \sigma) < \lambda$ for any cardinal $\sigma$ with $\kappa \leq \sigma < \lambda$. We need some preparation.

\medskip

\begin{Def}  For a cardinal $\theta \geq \lambda$ and $h : \theta \rightarrow P_\kappa(\lambda)$, we let $S^h_{\mu, \kappa, \lambda, \theta}$ denote the set of all $x \in P_\kappa(\theta)$ for which one can find $z \subseteq x$ such that 
\begin{itemize}
\item o.t.$(z) = \mu$ ;  
\item $h(\beta) \subset h(\alpha)$ whenever $\beta < \alpha$ are in $z$ ;    
\item  $x \cap \lambda \subseteq \bigcup_{\alpha \in z} h(\alpha)$.      
 \end{itemize}
\end{Def}

\begin{fact}  {\rm (\cite{Resemble})} 
\begin{enumerate}[\rm (i)] 
\item $S^h_{\mu, \kappa, \lambda, \theta} \subseteq X^h_{\mu, \kappa, \lambda, \theta}$.
\item Given $h : \theta \rightarrow P_\kappa(\lambda)$ and $x \in P_\kappa (\theta)$, the following are equivalent :
\begin{enumerate}[\rm (a)]
\item $x \in S^h_{\mu, \kappa, \lambda, \theta}$.
\item There is $\xi_i \in x$ for $i < \mu$ such that
\begin{itemize}
\item $h (\xi_j) \subseteq h (\xi_i)$ whenever $j < i < \mu$ ;
\item $x \cap \lambda \subseteq \bigcup_{i < \mu} h(\xi_i)$ ;
\item $(x \cap \lambda) \setminus h (\xi_i) \not= \emptyset$ for all $i < \mu$. 
\end{itemize}
\end{enumerate}
\end{enumerate}
 \end{fact}

\begin{fact}  {\rm (\cite{Resemble})} Suppose that $\cf (\lambda) = \mu$, $\langle \lambda_i : i < \mu \rangle$ is an increasing sequence of cardinals greater than or equal to $\kappa$ with supremum $\lambda$, and $s : \theta \rightarrow P_\kappa (\lambda)$ is such that $ran (s) \cap P_\kappa (\lambda_i) \in I^+_{\kappa, \lambda_i}$ for all $i < \mu$. Then $S^s_{\mu, \kappa, \lambda, \theta} \in (NG^\mu_{\kappa, \theta})^+$.
\end{fact}

\begin{Def}  We say that $s: \lambda\rightarrow P_\kappa(\lambda)$ is {\it easy} if $\bigcup_{\eta \in s (\xi)} s (\eta) \subseteq s (\xi)$ for all $\xi < \lambda$.
\end{Def}

\begin{fact}   {\rm(\cite{Secret})} Let $s : \lambda \rightarrow P_\kappa (\lambda)$ be easy and such that $X^s_{\mu, \kappa, \lambda, \lambda} \in I^+_{\kappa, \lambda}$. Then

\centerline{$N\mu$-$S_{\kappa,\lambda^+} = NS_{\kappa, \lambda^+} \vert (p^{- 1}_{\lambda^+} (X^s_{\mu, \kappa, \lambda, \lambda}) \cap \{ x \in P_\kappa (\lambda^+) : \cf (\sup x) = \mu\})$.}  
\end{fact}

\begin{Pro} Suppose that $\cf (\lambda)= \mu$, and $u (\kappa, \sigma) < \lambda$ for every cardinal $\sigma$ with $\kappa \leq \sigma < \lambda$. Then $Sup (N\mu$-$S_{\kappa,\lambda^+}) \subseteq NS_{\lambda^+} \vert E^{\lambda^+}_\mu$.
\end{Pro}

{\bf Proof.} Pick an increasing sequence $\langle \lambda_i : i < \mu \rangle$ of cardinals greater than or equal to $\kappa$ with supremum $\lambda$.

\medskip

{\bf Claim 1.}  There is an easy $s : \lambda \rightarrow \bigcup_{i < \mu} P_\kappa (\lambda_i)$ such that $ran (s) \cap P_\kappa (\lambda_i) \in I^+_{\kappa, \lambda_i}$ for all $i < \mu$.   

\medskip

{\bf Proof of Claim 1.}  Inductively define $\theta_i$ for $i < \mu$ as follows :
\begin{itemize}
\item $\theta_0 = 0$ ;
\item $\theta_{i + 1} = \max \{\theta_i^+, u (\kappa, \lambda_i)$ ;
\item $\theta_i = \sup \{ \theta_j : j < i \}$ in case $i$ is an infinite limit ordinal.
\end{itemize}
For $i < \mu$, pick $Z_i \in I^+_{\kappa, \lambda_i}$ with $\vert Z_i \vert = u (\kappa, \lambda_i)$, and select an onto function $r_i : \theta_{i + 1} \setminus \theta_i \rightarrow Z_i$. Set $r = \bigcup_{i < \mu} r_i$. By induction on $n < \omega$, define $a^n_\xi$ for $\xi < \lambda$ by
\begin{itemize}
\item $a^0_\xi = r (\xi)$ ;
\item $a^{n + 1}_\xi = a^n_\xi \cup (\bigcup_{\eta \in a^n_\xi} (\bigcup_{m \leq n} a^m_\eta))$.
\end{itemize}
Finally define $s : \lambda \rightarrow P_\kappa (\lambda)$ by $s (\xi) = \bigcup_{n < \omega} a^n_\xi$. It is simple to see that $s$ is easy. Now fix $i < \mu$.

Let us first show by induction that $A (n)$ holds for all $n < \omega$, where $A (n)$ asserts that $\bigcup_{m \leq n} a^m_\xi \subseteq \lambda_i$ for every $\xi \in \theta_{i + 1}$. It is immediate that $A (0)$ holds. Now suppose that $A (n)$ holds, and let $\xi \in \theta_{i + 1}$ and $\eta \in a^n_\xi$. Then $\eta$ belongs to $\lambda_i$, and hence to $\theta_{i + 1}$, so $\bigcup_{m \leq n} a^m_\eta \subseteq \lambda_i$. It  follows that $A (n + 1)$ holds.

Now let $a \in P_\kappa (\lambda_i)$. There must be $\xi < \theta_{i + 1}$ such that $a \subseteq r_i (\xi)$. Then

\centerline{$a \subseteq r (\xi) \subseteq s (\xi) =  \bigcup_{n < \omega} \bigcup_{m \leq n} a^m_\xi \subseteq \lambda_i$,}

which completes the proof of the claim. 

\medskip

By Claim 1 and Facts 4.5 (i) and 4.6, $X^s_{\mu, \kappa, \lambda, \lambda} \in I^+_{\kappa, \lambda}$. Now fix $S \in NS^+_{\lambda^+}$ with $S \subseteq E^{\lambda^+}_\mu$. Set $W = \{ x \in P_\kappa (\lambda^+) : Sup (x) \in S \}$ and $T = \{x \in W : x \cap \lambda \in X^s_{\mu, \kappa, \lambda, \lambda} \}$.

\medskip

{\bf Claim 2.}  $T \in NS^+_{\kappa, \lambda^+}$.   

\medskip

{\bf Proof of Claim 2.} Fix $F :\lambda^+ \times \lambda^+ \rightarrow \lambda^+$. We may find $\delta \in S \setminus \lambda$ with $F`` (\delta \times \delta) \subseteq \delta$. Pick an increasing sequence $\langle \delta_i : i < \mu \rangle$ of ordinals with supremum $\delta$. Now define inductively $x_i \in P_\kappa (\lambda^+)$ and $\alpha_i \in \lambda$ for $i < \mu$ so that
\begin{itemize}
\item $x_0 = \emptyset$ ;
\item $x_i \cap \lambda_i \subseteq s (\alpha_i) \subseteq \lambda_i$ ;
\item$x_{i + 1} = x_i \cup \{ \delta_i \} \cup F`` (x_i \times x_i) \cup \sup (x_i \cap \kappa) \cup \{ \alpha_i, \lambda_i \} \cup s (\alpha_i)$ ;
\item $x_i = \bigcup_{ j < i} x_j$ in case $i$ is an infinite limit ordinal.
\end{itemize}
Finally, put $x = \bigcup_{i < \mu} x_i$. Then clearly,
\begin{itemize}
\item $x \cap \kappa \in \kappa$ and $F`` (x \times x) \subseteq x$.
\item $Sup (x) \in S$.
\end{itemize}
Furthermore, $x \cap \lambda = \bigcup_{i < \mu} (x_i \cap \lambda_i) \subseteq \bigcup_{i < \mu} s (\alpha_i)$. Since $\{ \alpha_i : i < \mu \} \subseteq x$, and $s (\alpha_j) \subseteq x_{j + 1} \cap \lambda_j \subseteq x_i \cap \lambda_i \subseteq s (\alpha_i)$, it follows by Fact 4.5 that $x$ lies in $S^s_{\mu, \kappa, \lambda, \lambda}$, and hence in $X^s_{\mu, \kappa, \lambda, \lambda}$, which completes the proof of the claim. 

\medskip

It immediately follows from Claim 2 and Fact 4.8 that $W \in (N\mu$-$S_{\kappa,\lambda^+})^+$.
\hfill$\square$

\medskip

In Fact 4.3 and Proposition 4.9, the assumption is on the behaviour of the $u (\kappa, -)$ function. Let us recall a third result, with an assumption this time on the extent of Shelah's approachable ideal.

\medskip

\begin{Def}  Given a regular uncountable cardinal $\chi$, $I[\chi]$ denotes the collection of all $S \subseteq \chi$ for which there are a sequence $\langle u_\xi : \xi < \chi \rangle$ of subsets of $\chi$, and a closed unbounded subset $C$ of $\chi$ such that for any $\alpha \in S \cap C$,
\begin{itemize}
\item $\alpha$ is an infinite limit ordinal with $\cf(\alpha) < \alpha$ ;
\item there is a cofinal subset $B$ of $\alpha$ such that $o.t. (B) = \cf(\alpha)$ and $\{B \cap \eta : \eta < \alpha \} \subseteq \{ u_\xi : \xi < \alpha\}$.
\end{itemize}
\end{Def}

\medskip

The following is due to Shelah (for proofs and comments see \cite[Proposition 3.4, Theorem 3.18, Lemma 3.40, Corollary 4.6 and Proposition 4.11]{Eis}).

\medskip

\begin{fact} \begin{enumerate}[\rm (i)]
\item  Suppose that $\chi \notin I [\chi]$. Then $I[\chi]$ is a normal ideal on $\chi$.
\item  Let $\theta$ be a regular cardinal with $\theta^+ < \chi$. Then $NS_\chi^+ \cap P(E_\theta^\chi) \cap I[\chi] \not= \emptyset$.
\item Let $\theta$ be a regular cardinal such that $\rho^{<\theta} < \chi$ for every infinite cardinal $\rho < \chi$. Then $E_\theta^\chi \in I[\chi]$.
\item Suppose that $\chi$ is the successor of a regular cardinal $\rho$. Then for any regular cardinal $\theta < \rho$, $E_\theta^\chi \in I[\chi]$.
\item Suppose that $\chi$ is the successor of a cardinal $\rho$ such that  $\square_\rho^\ast$ holds. Then $\chi \in  I [\chi]$.
\end{enumerate}
\end{fact}

\medskip

The following is also due to Shelah (\cite{She93}, \cite{She1008}).

\medskip

\begin{fact} Given a regular cardinal $\theta < \chi$ and $S \subseteq E^\chi_\theta$, the following are equivalent :
\begin{enumerate}[\rm (i)]

\item  $S \in I[\chi]$.

\item  There are a sequence $\langle C_\xi : \xi < \chi \rangle$ and a closed unbounded subset $C$ of $\chi$ such that
\begin{itemize}
\item For each $\xi < \chi$, $C_\xi$ is a closed subset of $\xi$ ;
\item For any $\xi < \chi$, and any $\beta \in C_\xi$ with $\sup (C_\xi \cap \beta) < \beta$, $C_\beta = C_\xi \cap \beta$ ;
\item For any $\alpha \in S \cap C$, $o.t. (C_\alpha) = \theta$, and moreover $\sup C_\alpha = \alpha$.
\end{itemize}

\item There are a sequence $\langle z_\xi : \xi < \chi \rangle$ and a closed unbounded subset $C$ of $\chi$ such that
\begin{itemize}
\item For each $\xi < \chi$, $z_\xi \in P_\theta (\xi)$ ;
\item For any $\xi < \chi$, and any $\beta \in z_\xi$, $z_\beta = z_\xi \cap \beta$ ;
\item For any $\alpha \in S \cap C$, there is a cofinal subset $X$ of $\alpha$ of order-type $\theta$ such that $X \cap \xi = z_\xi$ for all $\xi \in X$. 
\end{itemize}

\end{enumerate}
\end{fact}

{\bf Proof.}

 \hskip0,4cm  (i) $\rightarrow$ (ii) : See \cite[Theorem 3.7]{Eis}.
  
\hskip0,2cm  (ii) $\rightarrow$ (i) : Given $\langle C_\xi : \xi < \chi \rangle$ and $C$ as in (ii), define $z_\xi$ for $\xi < \chi$ by : $z_\xi$ equals the collection of all $\gamma \in C_\xi$ with $\sup (C_\xi \cap \gamma) < \gamma$ if $\vert C_\xi \vert < \theta$, and $\emptyset$ otherwise.

\medskip

{\bf Claim 1.}  Let $\xi < \chi$, and $\beta \in z_\xi$. Then $z_\beta = z_\xi \cap \beta$.

\medskip

{\bf Proof of Claim 1.}  Since $\beta \in z_\xi$, $\sup (C_\xi \cap \beta) < \beta$, and therefore $C_\beta = C_\xi \cap \beta$. It follows that $\vert C_\beta \vert \leq \vert C_\xi \vert < \theta$. Hence \\$z_\beta =$ \\the collection of all $\gamma \in C_\beta$ with $\sup (C_\beta \cap \gamma) < \gamma =$ \\the collection of all $\gamma \in C_\xi \cap \beta$ with $\sup (C_\xi \cap \gamma) < \gamma = \\z_\xi \cap \beta$,\\ which completes the proof of the claim. 

\medskip

{\bf Claim 2.}  Let $\alpha \in S \cap C$.Then there is a cofinal subset $X$ of $\alpha$ of order-type $\theta$ such that $X \cap \xi = z_\xi$ for all $\xi \in X$. 

\medskip

{\bf Proof of Claim 2.}  Set $X = \{\xi \in C_\alpha : \sup (C_\alpha \cap \xi) < \xi\}$. Then clearly, $o.t. (X) = o.t. (C_\alpha) = \theta$, and moreover $\sup X = \alpha$. Now fix $\xi \in X$. Since $C_\xi = C_\alpha \cap \xi$, $\vert C_\xi \vert < \theta$, and consequently \\$z_\xi =$ \\the collection of all $\gamma \in C_\xi$ with $\sup (C_\xi \cap \gamma) < \gamma =$ \\ the collection of all $\gamma \in C_\alpha \cap \xi$ with $\sup (C_\alpha \cap \gamma) < \gamma =$\\ $X \cap \xi$, \\which completes the proof of the claim. 

\medskip

\hskip0,4cm  (iii) $\rightarrow$ (i) : Trivial.

 \hfill$\square$  

\begin{fact} {\rm(\cite{Ideals})} Suppose that $\lambda$ is regular. Then $Sup^{- 1} (S)$ lies in $(NG_{\kappa,\lambda}^\mu)^+$ (and hence by Fact 2.27 (ii) in $(N\mu$-$S_{\kappa,\lambda^+})^+$) for any $S \in NS_\lambda^+ \cap P (E_\mu^\lambda) \cap I [\lambda]$.  
\end{fact}
 
 \bigskip

\subsection{Towers}

\bigskip

Non-saturation of our ideals will be witnessed by the existence of long towers.

\medskip

\begin{Def} Given an ideal $K$ on a set $Q$, $Y\subseteq P (Q)$ and an ordinal $\delta$, a $(K, Y)${\it -tower of length} $\delta$ is a sequence $\langle A_\alpha : \alpha < \delta \rangle$ of subsets of $Q$ such that $(A_\alpha \setminus A_\beta, A_\beta \setminus A_\alpha) \in K^+ \times Y$ whenever $\alpha < \beta < \delta$.
\end{Def}

\begin{Obs} \begin{enumerate}[\rm (i)] 
\item Suppose that $\delta$ is a cardinal, and $\langle A_\alpha : \alpha < \delta \rangle$ is a $(K, Y)$-tower. Then $\langle A_\alpha \setminus A_{\alpha + 1} : \alpha < \delta \rangle$ witnesses that $K$ is not $Y$-$\delta$-saturated.
\item Suppose that for some $W \in K^+$, $\langle B_\alpha : \alpha < \delta \rangle$ is a $(K \vert W, K \vert W)$-tower. Then $\langle B_\alpha \cap W : \alpha < \delta \rangle$ is a $(K \vert W, K)$-tower.
\end{enumerate}
\end{Obs}

\begin{Obs}  Suppose that $\lambda$ is regular and for some ordinal $\delta$, there is a $(J, H)$-tower $\langle A_\alpha : \alpha < \delta \rangle$, where $J$ and $H$ are two ideals on $\lambda$. Let $K$ and $G$ be two ideals on $P_\kappa( \lambda)$ such that $Sup (K) \subseteq J$ and $H \subseteq Sup (G)$. Then $\langle Sup^{- 1} (A_\alpha) : \alpha < \delta \rangle$  is a $(K, G)$-tower.
\end{Obs}

{\bf Proof.} By Observation 4.2.
\hfill$\square$ 

\begin{Def}  Let $\sigma$ be a regular uncountable cardinal. Given $f, g \in {}^\sigma \sigma$, $f <_{I_\sigma} g$ means that $\vert \{\alpha < \sigma : f (\alpha) \geq g (\alpha \} \vert < \sigma$.

We let $\frak{b}_\sigma$ denote the least cardinality of any $F \subseteq {}^\sigma \sigma$ with the property that there is no $g \in {}^\sigma \sigma$ such that $f <_{I_\sigma} g$ for all $f \in F$. 

$Depth ({}^\sigma \sigma)$  denotes the least ordinal $\eta$ such that there is no increasing sequence $\langle f_i : i < \eta \rangle$ in $({}^\sigma \sigma, <_{I_\sigma})$.
 
We let $Depth ([\sigma]^\sigma, \searrow)$ denote the least ordinal $\eta$ such that there is no $(I_{\sigma}, I_{\sigma})$-tower of length $\eta$.
\end{Def}

\begin{fact}  \begin{enumerate}[\rm (i)] 
\item {\rm(\cite{She589})} $\frak{b}_\sigma < Depth ({}^\sigma \sigma)$.
\item {\rm (\cite{Towers})}  $Depth ({}^\sigma \sigma) \leq Depth ([\sigma]^\sigma, \searrow)$.
\end{enumerate}
\end{fact}

\medskip

We will produce $(N\mu$-$S_{\kappa,\lambda}, I_{\kappa, \lambda})$-towers of length $\frak{b}_{u (\kappa, \lambda)}$, and for $\tau < Depth ({}^{u (\kappa, \lambda)} u (\kappa, \lambda))$, $(N\mu$-$S_{\kappa,\lambda}, SNS_{\kappa, \lambda})$-towers of length $\tau$. In the case when $u (\kappa, \lambda) = \lambda = \cf (\lambda)$, we will start with towers on $\lambda$ (as opposed to $P_\kappa (\lambda)$) and then use this :


\begin{Lem}  Suppose that $\lambda$ is regular. Let $\tau$ be a cardinal, $K$ be a $\kappa$-complete ideal on $P_\kappa( \lambda)$, and $W \in NS^+_\lambda \cap P (E^\lambda_\mu)$ be such that $Sup (K) \subseteq NS_\lambda \vert W$. Then the following hold :
\begin{enumerate}[\rm (i)]
\item Suppose that there exists an $(NS_\lambda \vert W, NS_\lambda \vert W)$-tower of length $\tau$. Then there exists a $(K \vert Sup^{- 1} (W), SNS_{\kappa, \lambda})$-tower of length $\tau$. 
\item Suppose that there exists an $(NS_\lambda \vert W, I_\lambda)$-tower of length $\tau$. Then there exists a $(K \vert Sup^{- 1} (W), I_{\kappa, \lambda})$-tower of length $\tau$. 
\end{enumerate}
\end{Lem}

{\bf Proof.} Clearly, $Sup (K \vert Sup^{- 1} (W)) \subseteq NS_\lambda \vert W$. Now use Observations 4.2 (iii), 4.15 (ii) and 4.16 for (i), and Observations 4.2 (ii) and 4.16 for (ii).
\hfill$\square$

\medskip


By an argument of Gitik (see \cite{Towers}), if $J$ is a normal ideal on $\kappa$ that is not $\kappa^+$-saturated, then there exists a $(J, I_\kappa)$-tower of length $\kappa^+$. So a result formulated in terms of towers is not necessarily stronger. Gitik's result generalizes to $P_\kappa (\lambda)$ as follows.

\medskip

\begin{Obs} Let $J$ be a normal ideal on $P_\kappa (\lambda)$ that is not $\lambda^+$-saturated. Then there exists a $(J, I_{\kappa, \lambda})$-tower of length $\lambda^+$.
\end{Obs}

{\bf Proof.} Select $A_\alpha \in J^+$ for $\alpha < \lambda^+$ with the property that $A_\beta \cap A_\alpha \in J$ whenever $\beta < \alpha < \lambda^+$. For $\lambda \leq \alpha < \lambda^+$, pick a bijection $j_\alpha : \lambda \rightarrow \alpha$, and set $B_\alpha = \bigtriangleup_{i < \lambda} (P_\kappa (\lambda) \setminus A_{j_\alpha (i)})$. Notice that $B_\alpha \cap A_\beta \in I_{\kappa, \lambda}$ for all $\beta < \alpha$. 

\medskip

{\bf Claim 1.}  Let $\lambda \leq \alpha < \lambda^+$. Then $B_\alpha \in J^+$.

\medskip

{\bf Proof of Claim 1.}  It suffices to show that $A_\alpha \setminus B_\alpha \in J$. Suppose that this does not hold. Define $f : A_\alpha \setminus B_\alpha \rightarrow \lambda$ by $f (x)$ = the least $i \in x$ such that $x \in A_{j_\alpha (i)}$. We may find $i < \lambda$ and $H \in J^+ \cap P (A_\alpha \setminus B_\alpha)$ such that $f$ takes the constant value $i$ on $H$. But then $H \subseteq A_\beta \cap A_\alpha$, where $\beta = j_\alpha (i) < \alpha$. This contradiction completes the proof of the claim. 

\medskip

{\bf Claim 2.}  Let $\lambda \leq \beta < \alpha < \lambda^+$. Then $B_\beta \setminus B_\alpha \in J^+$.

\medskip

{\bf Proof of Claim 2.}  Clearly, $A_\beta \setminus (B_\beta \setminus B_\alpha) \subseteq (A_\beta \setminus B_\beta) \cup (A_\beta \cap B_\alpha)$, so $A_\beta \setminus (B_\beta \setminus B_\alpha) \in J$. The conclusion of the claim follows. 

\medskip

{\bf Claim 3.}  Let $\lambda \leq \beta < \alpha < \lambda^+$. Then $B_\alpha \setminus B_\beta \in J$.

\medskip

{\bf Proof of Claim 3.}  Supposing otherwise, define $g : B_\alpha \setminus B_\beta \rightarrow \lambda$ by $g (x)$ = the least $i \in x$ such that $x \in A_{j_\beta (i)}$. We may find $i < \lambda$ and $G \in J^+ \cap P (B_\alpha \setminus B_\beta)$ such that $g$ is constantly $i$ on $G$. There must be $k < \lambda$ such that $j_\alpha (k) = j_\beta (i)$. But then $G \subseteq \{x \in P_\kappa (\lambda) : k \notin x \}$. This contradiction completes the proof of the claim. 

\medskip

For $\lambda \leq \beta < \alpha < \lambda^+$, pick $C_{\beta\alpha} \in J^\ast$ with $(B_\alpha \setminus B_\beta) \cap C_{\beta\alpha} = \emptyset$. By induction on $n$, define $D^n_\alpha$ for $n < \omega$ and $\lambda \leq \alpha < \lambda^+$ so that
\begin{itemize}
\item $D^0_\alpha \setminus C_{\beta\alpha} \in I_{\kappa, \lambda}$ whenever $\lambda \leq \beta < \alpha$ ;
\item  $D^{n + 1}_\alpha \setminus D^n_\beta \in I_{\kappa, \lambda}$ whenever $\lambda \leq \beta < \alpha$.
\end{itemize}
For  $\lambda \leq \alpha < \lambda^+$, put $D_\alpha = \bigcap_{n < \omega} D^n_\alpha$. Notice that $D_\alpha \in J^\ast$.

\medskip

{\bf Claim 4.}  Let $\lambda \leq \beta < \alpha < \lambda^+$. Then $(B_\beta \cap D_\beta) \setminus (B_\alpha \cap D_\alpha) \in J^+$.

\medskip

{\bf Proof of Claim 4.}  Clearly, $(B_\beta \setminus B_\alpha) \cap D_\beta \subseteq (B_\beta \cap D_\beta) \setminus (B_\alpha \cap D_\alpha)$. By Claim 2, the conclusion of the claim follows. 

\medskip

{\bf Claim 5.}  Let $\lambda \leq \beta < \alpha < \lambda^+$. Then $(B_\alpha \cap D_\alpha) \setminus (B_\beta \cap D_\beta) \in I_{\kappa, \lambda}$.

\medskip

{\bf Proof of Claim 5.}  It is simple to see that $(B_\alpha \cap D_\alpha) \setminus (B_\beta \cap D_\beta) = S \cup \bigcup_{n < \omega} T_n$, where $S = (B_\alpha \cap D_\alpha) \setminus B_\beta$ and $T_n = (B_\alpha \cap D_\alpha) \setminus D^n_\beta$. Now, $S \in I_{\kappa, \lambda}$ since $S \subseteq (B_\alpha \setminus B_\beta) \cap D^0_\alpha$. Furthermore for each $n < \omega$, $T_n \in I_{\kappa, \lambda}$ since $T_n \subseteq D^{n + 1}_\alpha \setminus D^n_\beta$. The conclusion of the claim follows.
\hfill$\square$

\bigskip

\subsection{Case when $\lambda$ is the successor of a cardinal $\nu$ with $\cf (\nu) \not= \mu$}
 
\bigskip
 
\begin{fact}  {\rm(\cite{Towers})} Suppose that $\lambda = \nu^+$, where $\cf (\nu ) \not= \mu$. Then the following hold :
\begin{enumerate}[\rm (i)]
\item For any regular cardinal $\tau < Depth ({}^\lambda \lambda)$, and any $\lambda$-complete ideal $J$ on $\lambda$ extending $NS_\lambda \vert E^\lambda_\mu$, there exists a $(J, J)$-tower of length $\tau$. 
\item For any $S \in NS^+_\lambda \cap P (E^\lambda_\mu)$, there exists an $(NS_\lambda \vert S, I_\lambda)$-tower of length $\frak{b}_\lambda$. 
\end{enumerate}
\end{fact}

\begin{Obs}  Suppose that $\lambda = \nu^+ = 2^\nu$, where $\cf (\nu ) \not= \mu$. Then for any $\eta < Depth ([\lambda]^\lambda, \searrow)$, and any  $\lambda$-complete ideal $J$ on $\lambda$ extending $NS_\lambda \vert E^\lambda_\mu$, there exists a $(J, I_\lambda)$-tower of length $\eta$. 
\end{Obs}

{\bf Proof.} Let $0 < \eta < Depth ([\lambda]^\lambda, \searrow)$, and let $J$ be a complete ideal $J$ on $\lambda$ with $NS_\lambda \vert E^\lambda_\mu \subseteq J$. By a result of Shelah \cite{She10}, $\diamondsuit_\lambda [J]$ holds. By a result of \cite{Towers}, the desired conclusion follows.
\hfill$\square$ 

\begin{Rmk} If, as in the assumption of Observation 4.22, $2^{< \lambda} = \lambda = \cf (\lambda)$, then there exists an almost disjoint family of $2^\lambda$ subsets of $\lambda$. Does it follow that $2^\lambda < Depth ([\lambda]^\lambda, \searrow)$ ?
\end{Rmk}

\begin{Pro}  Suppose that $\lambda = \nu^+$, where $\cf (\nu ) \not= \mu$, and let $W \in NS^+_\lambda \cap P (E_\mu^\lambda) \cap I [\lambda]$. Then the following hold :
\begin{enumerate}[\rm (i)]
\item For any regular cardinal $\tau < Depth ({}^\lambda \lambda)$, there exists an $(NG_{\kappa,\lambda}^\mu \vert Sup^{- 1} (W), SNS_{\kappa, \lambda})$-tower of length $\tau$. 
\item There exists an $(NG_{\kappa,\lambda}^\mu \vert Sup^{- 1} (W), I_{\kappa, \lambda})$-tower of length $\frak{b}_\lambda$. 
\end{enumerate}
\end{Pro}

{\bf Proof.} By Lemma 4.19 and Facts 4.13 and 4.21.
\hfill$\square$ 

\medskip

\begin{Rmk} \begin{enumerate}[\rm (i)]
\item By Fact 4.11 (ii), it is always possible to find such a $W$.
\item If $\nu$ is regular, then by Fact 4.11 (iv), we can take for $W$ any stationary subset of $E_\mu^\lambda$.
\item Notice that by Fact 2.27 (ii), any $(NG_{\kappa,\lambda}^\mu \vert Sup^{- 1} (W), SNS_{\kappa, \lambda})$-tower (respectively, $(NG_{\kappa,\lambda}^\mu \vert Sup^{- 1} (W), I_{\kappa, \lambda})$-tower) of length $\tau$ is an \\ $(N\mu$-$S_{\kappa,\lambda} \vert Sup^{- 1} (W), SNS_{\kappa, \lambda})$-tower (respectively, $(N\mu$-$S_{\kappa,\lambda} \vert Sup^{- 1} (W), I_{\kappa, \lambda})$-tower) of length $\tau$.
\end{enumerate}
\end{Rmk}

\medskip

The condition that $W \in I [\lambda]$ can be dropped if we strengthen the assumption.

\medskip

\begin{Pro}  Suppose that $\lambda = \nu^+$, where $u (\kappa, \nu) = \nu$. Then for any $W \in NS^+_\lambda \vert E^\lambda_\mu$, the following hold :
\begin{enumerate}[\rm (i)]
\item For any regular cardinal $\tau < Depth ({}^\lambda \lambda)$, there exists a $(N\mu$-$S_{\kappa,\lambda}  \vert Sup^{- 1} (W), SNS_{\kappa, \lambda})$-tower of length $\tau$. 
\item There exists a $(N\mu$-$S_{\kappa,\lambda}  \vert Sup^{- 1} (W), I_{\kappa, \lambda})$-tower of length $\frak{b}_\lambda$. 
\end{enumerate}
\end{Pro}

{\bf Proof.} By Fact 1.9, $\cf (\nu) \geq \kappa$. Now use Lemma 4.19 and Facts 4.3 and 4.21.
\hfill$\square$

\bigskip

\subsection{Case when $\lambda$ is either weakly inaccessible, or the successor of a cardinal $\nu$ with $\cf (\nu) = \mu$}
 
\bigskip
 
\begin{Def}  Let $\sigma$ be a regular uncountable cardinal, and $W$ be a stationary subset of $\sigma$. For $\gamma \in E^\sigma_{\geq \omega_1}$, $W$ {\it reflects at} $\gamma$ if $W \cap \gamma$ is stationary in $\gamma$.

Given a stationary subset $T$ of $E^\sigma_{\geq \omega_1}$, $W$ {\it reflects fully in} $T$ if there is a club subset $C$ of $\sigma$ such that $W$ reflects at every $\gamma \in G \cap T$.
 \end{Def}

\medskip

We will rely on the following \say{towerized} version of a celebrated result of Gitik and Shelah \cite{GS}. 

\medskip

\begin{fact} {\rm(\cite{Towers})}  Suppose that $\tau < \rho < \sigma$ are three regular cardinals, and $W$ is a stationary subset of $E^\sigma_\tau$ that reflects fully in $E^\sigma_{\geq \rho}$. Suppose further that either $\tau > \omega$, or $\rho \geq \omega_2$. Then there is an $(NS_\sigma \vert W, I_\sigma)$-tower of length ${\frak b}_\sigma$. 
\end{fact}

\medskip

The following, which is due to Shelah, is immediate from Theorem 0.1 in \cite{She1008} and Fact 4.12.

\medskip

\begin{fact} Let $\tau, \rho, \sigma$ be three regular uncountable cardinals with $\tau^{++} < \rho < \sigma$. Then there is a $W$ in $I [\sigma] \cap P (E^\sigma_\tau)$ that reflects at every $\gamma \in E^\sigma_\rho$. 
\end{fact}

\begin{Cor} Let $\tau, \rho, \sigma$ be three regular uncountable cardinals with $\tau^{++} < \rho < \sigma$. Then there is a stationary $W$ in $I [\sigma] \cap P (E^\sigma_\tau)$ that reflects fully in $E^\sigma_{\geq \rho}$. 
\end{Cor}

{\bf Proof.} By Fact 4.29, there is a $W$ in $I [\sigma] \cap P (E^\sigma_\tau)$ that reflects at every $\gamma \in E^\sigma_\rho$. It is simple to see that $W$ is stationary. Furthermore $W$ reflects fully in $E^\sigma_\rho$, and hence in $E^\sigma_\chi$ for every regular cardinal $\chi$ with $\rho \leq \chi < \sigma$.
\hfill$\square$

\begin{Pro} Let $\tau , \sigma$ be two regular uncountable cardinals with $\tau^{+ 3} < \sigma$. Then there is an $(NS_\sigma \vert W, I_\sigma)$-tower of length ${\frak b}_\sigma$ for some stationary $W$ in $I [\sigma] \cap P (E^\sigma_\tau)$. 
\end{Pro}

{\bf Proof.} By Fact 4.28 and Corollary 4.30.
\hfill$\square$

\begin{Pro}  Suppose that either $\lambda$ is weakly inaccessible, or $\lambda = \nu^+$, where $\cf (\nu) = \mu$. Then there exists an $(NG^\mu_{\kappa,\lambda}, I_{\kappa, \lambda})$-tower of length $\frak{b}_\lambda$. 
\end{Pro}

{\bf Proof.} Case when  $\mu > \omega$. Then Proposition 4.31 tells us that there is an $(NS_\lambda \vert W, I_\lambda)$-tower of length ${\frak b}_\lambda$ for some stationary $W$ in $I [\lambda] \cap P (E^\lambda_\mu)$. By Fact 4.13 and Lemma 4.19, the existence of an $(NG^\mu_{\kappa,\lambda} \vert Sup^{- 1} (W), I_{\kappa, \lambda})$-tower of length $\frak{b}_\lambda$ follows. 

\medskip

Case when  $\mu = \omega$. Then clearly, $E^\lambda_\mu$ reflects at every $\gamma \in E\lambda_{\geq \omega_2}$, so by Fact 4.28, there must be an $(NS_\lambda \vert E^\lambda_\mu, I_\lambda)$-tower of length $\frak{b}_\lambda$. Furthermore by Fact 4.11 (ii), $E^\lambda_\mu$ lies in $I [\lambda]$. By Fact 4.13 and Lemma 4.19, it follows that there exists an $(NG^\mu_{\kappa,\lambda}, I_{\kappa, \lambda})$-tower of length $\frak{b}_\lambda$. 
\hfill$\square$


\bigskip



\subsection{Case when $\cf (\lambda) \in \kappa \setminus \{\mu\}$}
 
 \bigskip
 
Throughout this subsection we let $\pi$ denote a regular cardinal greater than $\lambda$.

\medskip
 
\begin{Obs}  Suppose that for some ordinal $\delta$, there is a $(J, H)$-tower $\langle A_\alpha : \alpha < \delta \rangle$, where $J$ and $H$ are two ideals on $P_\kappa (\pi)$. Let $K$ and $G$ be two ideals on $P_\kappa( \lambda)$ such that $\chi (K) \subseteq J$ and $H \subseteq \chi (G)$, where $\chi : P_\kappa (\lambda) \rightarrow P_\kappa (\pi)$. Then $\langle \chi^{- 1} (A_\alpha) : \alpha < \delta \rangle$  is a $(K, G)$-tower.
\end{Obs}


\begin{Pro}  \begin{enumerate}[\rm (i)]
\item Suppose that there exists a $(\mu, \kappa, \lambda, \pi)$-shuttle $\chi$. Then there exists an $(N\mu$-$S_{\kappa,\lambda}, I_{\kappa, \lambda})$-tower of length $\frak{b}_\pi$. 
\item Suppose that $\pi = \nu^+$, where $\cf (\nu) \not= \mu$, and there exists a $(\mu, \kappa, \lambda, \pi)$-shuttle $\chi$. Then for any regular cardinal $\tau < Depth ({}^\pi \pi)$, there exists an $(N\mu$-$S_{\kappa,\lambda}, N\mu$-$S_{\kappa,\lambda})$-tower of length $\tau$. 
\item Suppose that $\pi = \nu^+$, where $\cf (\nu) \not= \mu$, and ${\mathcal A}_{\kappa,\lambda} (\mu,\pi)$ holds. Then for any regular cardinal $\tau < Depth ({}^\pi \pi)$, there exists an $(NG^\mu_{\kappa,\lambda}, NS_{\mu, \kappa,\lambda})$-tower of length $\tau$.
\end{enumerate}
\end{Pro}

{\bf Proof.} 

(i) : By Fact 2.27 (ii), Remark 4.25 (i) and Propositions 4.24 and 4.32, there exists an $(N\mu$-$S_{\kappa, \pi}, I_{\kappa, \pi})$-tower of length $\frak{b}_\pi$. By Proposition 3.7 and Observations 3.3 and 4.33, the existence of an $(N\mu$-$S_{\kappa,\lambda}, I_{\kappa, \lambda})$-tower of length $\frak{b}_\pi$ follows. 

\medskip

(ii) : Let $\tau < Depth ({}^\pi \pi)$ be a regular cardinal. Then by Fact 2.27 (ii), Remark 4.25 (i) and Proposition 4.24, there exists an $(N\mu$-$S_{\kappa, \pi}, SNS_{\kappa, \pi})$-tower of length $\tau$. Notice that by Proposition 3.7, $SNS_{\kappa, \pi} \subseteq N\mu$-$S_{\kappa,\pi} = \chi (N\mu$-$S_{\kappa,\lambda})$. Hence by Observation 4.33, there must be an $(N\mu$-$S_{\kappa,\lambda},N\mu$-$S_{\kappa,\lambda})$-tower of length $\tau$. 


 \medskip
 

(iii) : Proceed as in the proof of (ii), using Propositions 3.27, 3.39 and 3.40.
\hfill$\square$

\bigskip

\subsection{Summary}
 
\bigskip

Theorem 1.2 will be derived from the following. 

\medskip

\begin{Pro}  \begin{enumerate}[\rm (i)]
\item Assume that $\lambda$ is regular. Then there exists an $(NG^\mu_{\kappa,\lambda}, I_{\kappa, \lambda})$-tower of length $\frak{b}_\lambda$. 
\item Assume that one of the following holds :
\begin{enumerate}[\rm (i)]
\item $\cf (\lambda) < \mu$ and ${\mathcal A}_{\kappa,\lambda} (\mu, \lambda^+)$.
\item $\mu < \cf (\lambda) < \kappa$ and for some infinite set $A$ of regular infinite cardinals with $\vert A \vert < \min (\kappa, \min (A))$ and $\sup A = \lambda$, some $\mu^+$-complete ideal $I$ on $A$, and some $<_I$-increasing, cofinal sequence $\vec{g} = \langle g_\xi : \xi < \lambda^+ \rangle$ in $(\prod A, <_I)$, there is a closed unbounded subset $C$ of $\lambda^+$ such that every $\delta$ in $C \cap E^{\lambda^+}_\kappa$  is a good point for $\vec g$. 
\end{enumerate}
Then there exists an $(N\mu$-$S_{\kappa,\lambda}, I_{\kappa, \lambda})$-tower of length $\frak{b}_{\lambda^+}$. 
\end{enumerate}
\end{Pro}

{\bf Proof.} (i) : By Remark 4.25 (i) and Propositions 4.24 and 4.32.

\medskip

(ii) : By Propositions 3.52 and 4.34 (i).
\hfill$\square$

\bigskip

\section{Paradise not found}
 
\bigskip

In this section, we discuss the remaining cases.

\bigskip

\subsection{The case $\kappa \leq \cf (\lambda) < \lambda$}

\bigskip

Assuming that SSH holds and $\kappa \leq \cf (\lambda) < \lambda$, we have that $u (\kappa, \lambda) = \lambda$, so $N\mu$-$S_{\kappa,\lambda}$ and $N\mu$-$S_{\kappa, u (\kappa, \lambda)}$ are (equal and hence trivially) isomorphic. As for towers, it is not clear which length we should go for. Bounding numbers do not seem relevant, since $u (\kappa, \lambda)$ is not regular. What should be expected is, as in \cite{ForMag2}, the use of scales on $\lambda$, which is necessarily tricky since members of such scales do not have size less than $\kappa$. Our proof being based on mutual stationarity and Namba combinatorics on trees of countable height, it only deals with the case $\mu = \omega$.

\medskip

\begin{fact}  {\rm(\cite{Game})} There is $y : \lambda^{< \kappa} \rightarrow P_\kappa (\lambda)$ such that (a) for any $\beta \in \lambda$, $\beta \in y (\beta)$, and (b) $NG^\omega_{\kappa, \lambda} = {\overline y} (NS^+_{\omega_1, \lambda^{< \kappa}})$, where ${\overline y} : P_{\omega_1} (\lambda^{< \kappa}) \rightarrow P_\kappa (\lambda)$ is defined by ${\overline y} (a) = \bigcup_ {\beta \in a} y (\beta)$.
\end{fact}

\begin{Def}  Let $\sigma < \tau$ be two infinite cardinals, and $I$ be an ideal on $\sigma$. We put $\PP_I^\ast(\tau)$ denote the set of all cardinals $\theta$ for which one can find a sequence $\langle \tau_i : i < \sigma \rangle$ of regular infinite cardinals less than $\tau$ with supremum $\tau$ such that $\{i < \sigma : \tau_i  \leq \xi\} \in I$ for all $\xi < \tau$, and $\tcf(\prod_{i < \sigma} \tau_i, <_I) = \theta$.

\end{Def}
 

\medskip

The following is a version with towers of a result in \cite{Nonsat2}. Under SSH, all it gives is a tower of length $\lambda^+$.

\medskip

\begin{Pro} Suppose that $\kappa \leq \cf (\theta) < \theta \leq \lambda$. Then for any $\pi \in \PP_I^\ast(\theta)$, where $I = I_{\cf (\theta)}$, there is an $(NG^\omega_{\kappa,\lambda}, K)$-tower of length $\pi$, where $K$ is the $\cf (\theta)$-complete ideal on $P_\kappa (\lambda)$ defined by 

\centerline{$K = \bigcup_{r < \cf (\theta)} P (\{x \in P_\kappa (\lambda) : x \cap \cf (\theta) \subseteq r \} )$.} 
\end{Pro}

{\bf Proof.}   Put $\tau = \cf (\theta)$ and $\nu = \lambda^{< \kappa}$. For $\delta \in E^\tau_\omega$, pick an order-type $\omega$, cofinal subset $e_\delta$ of $\delta$. By a result of \cite{SSH}, we may find an increasing (and hence one-one) sequence $\langle \theta_i : i < \tau \rangle$ of regular infinite cardinals less than $\theta$ with supremum $\theta$ such that $\{i < \tau : \theta_i  \leq \xi\} \in I$ for all $\xi < \theta$, and $\tcf(\prod_{i < \tau} \theta_i, <_I ) = \pi$. For $i < \tau$, select a sequence $\langle S^{i}_\alpha : \alpha < \theta_i \rangle$ of disjoint stationary subsets of $E^{\theta_i}_\omega$. For $h \in \prod_{i < \tau} \theta_i$, let $Q_h$ denote the set of all $a \in P_{\omega_1} (\nu)$ such that (a) $\sup (a \cap \tau) \notin a$ and $e_{\sup (a \cap \tau)} \subseteq a$, and (b) $\sup (a \cap \theta_i) \in S^{i}_{h (i)}$ for all $i \in a \cap \tau$.

\medskip

{\bf Claim 1. }  $Q_h \in  NS^+_{\omega_1, \nu}$.  

\medskip

{\bf Proof of Claim 1.}   Let $\langle \sigma_j : j < \delta \rangle$ be the increasing enumeration of the set $\{ \tau \} \cup \{ \theta_i : i < \tau \}$. Define $k : \delta \rightarrow \tau + 1$ by : $k (j)$ equals the unique $i < \tau$ such that $\sigma_j = \theta_i$ if $\sigma_j \in \{ \theta_i : i < \tau \}$, and $\tau$ otherwise. Let $\tau = h (r)$. We let $D$ demote the set of all $a \in P_{\omega_1} (\nu)$ such that (a) $r \in a$, and (b) $h`` (a \cap \delta) = a \cap \tau$. Notice that $D \in NS^\ast_{\omega_1, \nu}$. Define $s \in \prod_{j < \delta} P (\sigma_j)$ by : $s (j)$ equals $S^{\sigma_j}_{h (k (j))}$ if $\sigma_j \in \{ \theta_i : i < \tau \}$, and $E^\tau_\omega$ otherwise.  Let $A$ be the set of all $a \in P_{\omega_1} (\nu)$ such that (a) $\sup (a \cap \tau) \notin a$ and $e_{\sup (a \cap \tau)} \subseteq a$, and (b) $\sup (a \cap \sigma_j) \in s (j)$ for all $j \in a \cap \delta$. By Proposition 2.1 of \cite{Nonsat2}, $A \in  NS^+_{\omega_1, \nu}$. We claim that $D \cap A \subseteq Q_h$. Thus let $a \in D \cap A$ and $i \in a \cap \tau$. Let $\theta_i = \sigma_j$. Then $\sup (a \cap \sigma_j) \in s(j)$, since $j \in a \cap \delta$, which is the same as saying that  $\sup (a \cap \theta_i) \in S^{i}_{h (i)}$. This completes the proof of the claim.

\medskip

Let $\langle f_\alpha : \alpha < \pi \rangle$ be an increasing, cofinal sequence in $(\prod_{i < \tau} \theta_i , <_I)$. Put $c = \{ i < \tau : \theta_i < \kappa \}$. Note that $c \in I$. For $\alpha < \pi$, let $W_\alpha$ be the set of all $a \in P_{\omega_1} (\nu)$ such that  (a) $\sup (a \cap \tau) \notin a$ and $e_{\sup (a \cap \tau)} \subseteq a$, and (b) there is $\xi \in a \cap \tau$ such that $\sup (a \cap \theta_i) \in \bigcup \{S^{i}_\gamma : f_\alpha (i) \leq \gamma < \theta_i \}$ for all $i \in (a \cap \tau) \setminus (c \cup \xi)$.

Let $y$ and ${\overline y}$ be as in the statement of Fact 5.1.  We let $C$ denote the set of all infinite $a \in P_{\omega_1} (\nu)$ such that (a) $\sup c < \sup (a \cap \tau)$, (b) for any $\eta \in a$, there is $\beta \in a \cap \tau$ with $\sup (y (\eta) \cap \tau) < \beta$, and (c)  for any $\eta \in a$ and any $i \in a \cap \tau$, there is $\beta \in a \cap \theta_i$ with $\sup (y (\eta) \cap \theta_i) < \beta$. Then clearly, $C$ is a closed unbounded subset of $P_{\omega_1} (\nu)$. Notice that for each $a \in C$, $\sup ({\overline y} (a) \cap \tau) = \sup (a \cap \tau)$, and moreover  $\sup ({\overline y} (a) \cap \theta_i) = \sup (a \cap \theta_i)$ for all $i \in a \cap \tau$. For $\alpha < \pi$, set $X_\alpha = {\overline y}`` (C \cap W_\alpha)$.

\medskip

{\bf Claim 2. }  Let $\alpha < \beta < \pi$. Then $X_\alpha \setminus X_\beta \in (NG^\omega_{\kappa, \lambda})^+$.

\medskip

{\bf Proof of Claim 2.} Let $r < \tau$ be such that $\{ i < \tau : f_\alpha (i) \geq f_\beta (i) \} \subseteq r$. Put $D = \{a \in P_{\omega_1} (\nu) : \sup (a \cap \tau) > r \}$, $A = D \cap C \cap Q_{f (\alpha)})$ and $Y = {\overline y}`` (A)$. Clearly, $Y \subseteq X_\alpha$, and moreover by Claim 1, $Y \in  (NG^\omega_{\kappa, \lambda})^+$. We claim that $Y \cap X_\beta = \emptyset$. Suppose otherwise, and pick $x \in Y \cap X_\beta$. We may find $a \in A$ and $b \in C \cap W_\beta$ such that ${\overline y} (a) = x = {\overline y} (b)$, and $\xi \in b \cap \tau$ such that $\sup (b \cap \theta_i) \in \bigcup \{S^{i}_\gamma : f_\beta (i) \leq \gamma < \theta_i \}$ for all $i \in (b \cap \tau) \setminus (c \cup \xi)$. Notice that $a \cup b \subseteq x$. Put $\delta = \sup (x \cap \tau)$. Then $\sup (a \cap \tau) = \delta = \sup (b \cap \tau)$, and moreover $e_\delta \subseteq x$. Now pick $i \in e_\delta$ with $i \notin c \cup r \cup \xi$. Then $\sup (a \cap \theta_i) = \sup (x \cap \theta_i) = \sup (b \cap \theta_i)$. However, $\sup (a \cap \theta_i) \in S^{i}_{f_\alpha (i)}$, whereas $\sup (b \cap \theta_i) \in \bigcup \{S^{i}_\gamma : f_\beta (i) \leq \gamma < \theta_i \}$. This contradiction completes the proof of the claim.

\medskip

{\bf Claim 3. }  Let $\alpha < \beta < \pi$, and let $r < \tau$ be such that $\{ i < \tau : f_\alpha (i) \geq f_\beta (i) \} \subseteq r$. Then $X_\beta \setminus X_\alpha \subseteq \{x \in P_\kappa (\lambda) : x \cap \tau \subseteq r \}$.

\medskip

{\bf Proof of Claim 3.}   Fix $x \in X_\beta \setminus \{x \in P_\kappa (\lambda) : x \cap \tau \subseteq r \}$. There must be $a \in C \cap W_\beta$ such that ${\overline y} (a) = x$. Notice that $\sup (a \cap \tau) = \sup (x \cap \tau)$. We may find $\xi \in a \cap \tau$ such that $\sup (a \cap \theta_i) \in \bigcup \{S^{i}_\eta : f_\beta (i) \leq \eta < \theta_i \}$ for all $i \in (a \cap \tau) \setminus (c \cup \xi)$. Pick $\zeta \in (a \cap \tau) \setminus r$ with $\zeta \geq \xi$. Then clearly,  $\sup (a \cap \theta_i) \in \bigcup \{S^{i}_\eta : f_\alpha (i) \leq \eta < \theta_i \}$ for all $i \in (a \cap \tau) \setminus (c \cup \zeta)$. Thus $a \in W_\alpha$. It follows that $x \in X_\alpha$, which completes the proof of the claim.

\medskip

By Claims 2 and 3, $\langle X_\alpha : \alpha < \pi \rangle$ is an $(NG^\omega_{\kappa,\lambda}, K)$-tower.
\hfill$\square$

\bigskip

\subsection{The case $\cf (\lambda) = \mu$}
 
\bigskip

We start with the following observation.

\medskip

\begin{Obs}  Let $S \in NS^+_{\kappa, \lambda}$. Then there is no normal, fine ideal $K$ on $P_\kappa (2^\lambda)$ such that $NS_{\kappa, \lambda} \vert S$ and $K$ are isomorphic.
\end{Obs}

{\bf Proof.} Suppose that $NS_{\kappa, \lambda} \vert S$ is isomorphic to some normal, fine ideal $K$ on $P_\kappa (2^\lambda)$. Then by Observation 2.8 and Fact 2.13, $\cof (K) = \cof (NS_{\kappa, \lambda} \vert S) \leq 2^\lambda$. This is a contradiction, since by a result of \cite{MPS1}, every normal, fine ideal on $P_\kappa (2^\lambda)$ has cofinality greater than $2^\lambda$.
\hfill$\square$

\begin{Def}  We let $\Phi (\mu, \kappa, \lambda)$ denote the least cardinal $\theta \geq \lambda$ for which there is $h : \theta \rightarrow P_\kappa(\lambda)$ such that $X^h_{\mu, \kappa, \lambda, \theta} \in I^+_{\kappa, \theta}$.
\end{Def}

\begin{fact} \begin{enumerate}[\rm (i)]
\item {\rm (\cite{Secret})} Suppose that $\mu = \omega$ and either $\kappa = \omega_1$, or $\lambda < \kappa^{+ \omega_1}$. Then 

\centerline{$NS_{\kappa,\lambda} \vert E^{\kappa, \lambda}_\mu = NS_{\mu, \kappa,\lambda} \vert E^{\kappa, \lambda}_\mu = N\mu$-$S_{\kappa,\lambda} = NG_{\kappa,\lambda}^\mu$.}

\item {\rm (\cite{Secret})} Assume that $\Phi (\mu, \kappa, \lambda) = \lambda$. Then there is $S \in NS^+_{\kappa, \lambda}$ such that $N\mu$-$S_{\kappa,\lambda} = NS_{\kappa,\lambda} \vert S = SNS_{\kappa,\lambda} \vert S$. 
\item {\rm (\cite{Resemble})} Assuming SSH, $\Phi (\mu, \kappa, \lambda)$ equals $\lambda^+$ if $\cf (\lambda) \in \kappa \setminus \{ \mu \}$, and $\lambda$ otherwise.
\end{enumerate}
\end{fact} 

\medskip

Thus, assuming GCH and $\cf (\lambda) = \mu$, $N\mu$-$S_{\kappa,\lambda}$ and $N\mu$-$S_{\kappa, u (\kappa, \lambda)}$ are not isomorphic, and in fact no restriction of $N\mu$-$S_{\kappa,\lambda}$ is isomorphic to a normal ideal on $P_\kappa (u (\kappa, \lambda))$. Under the same assumptions, we can show that if there is no inner model with large large cardinals, then some restriction of $NG_{\kappa, u (\kappa, \lambda)}^\mu$ is isomorphic to a normal extension of $NG_{\kappa,\lambda}^\mu$.

\medskip

Throughout the remainder of this subsection we let $\pi$ denote a regular cardinal greater than $\lambda$.

\medskip

\begin{Def} Let $\tau$ be a cardinal with $2 \leq \tau \leq \kappa$. A ${\mathcal B}_{\kappa,\lambda} (\tau,\pi)${\it -sequence} is a $(\tau,\lambda,\pi)$-sequence $\vec y = \langle y_\beta : \beta < \pi \rangle$ with the property that for each nonempty $e$ in $P_{\kappa^+} (\pi)$, there is a $< \kappa$-to one $g \in \prod_{\beta \in e} y_\beta$.

${\mathcal B}_{\tau,\lambda} (\kappa,\pi)$ asserts the existence of a ${\mathcal B}_{\kappa,\lambda} (\tau,\pi)$-sequence.
\end{Def}

\medskip

It is simple to see that any ${\mathcal B}_{\kappa,\lambda} (\tau,\pi)$-sequence is an ${\mathcal A}_{\kappa,\lambda} (\tau,\pi)$-sequence.

\medskip

There is in \cite{Norm} a result that is similar to Observation 3.51. From its stronger assumption (that every $\delta$ in $C \cap E^\pi_\kappa$  is (not just a good point, but) a \underline{remarkably} good point for $\vec g$), it is derived that ${\mathcal B}_{\sigma^+, \lambda} (\kappa,\pi)$ holds.

\medskip

\begin{fact} {\rm (\cite{Ideals})}
\begin{enumerate}[\rm (i)]
\item Suppose that $\lambda^{< \mu} < \pi$, and $\vec y = \langle y_\beta : \beta < \pi \rangle$ is an ${\mathcal A}_{\kappa,\lambda} (\mu^+,\pi)$-sequence. Then ($\{x \in P_\kappa (\pi) : \sup x \in f_{\vec y} (x \cap \lambda)\}$, and hence) $P_\kappa (\pi) \setminus A_\kappa (\vec y)$ lies in $(NG_{\kappa,\pi}^\mu)^+$.
\item Let $\vec y = \langle y_\beta : \beta < \pi \rangle$ be a ${\mathcal B}_{\kappa,\lambda} (\kappa,\pi)$-sequence. Then $A_\kappa (\vec y) \in (NG_{\kappa,\pi}^\mu)^+$.
\end{enumerate}
\end{fact} 

\begin{Obs} Let $\vec y = \langle y_\beta : \beta < \pi \rangle$ be a ${\mathcal B}_{\kappa,\lambda} (\kappa,\pi)$-sequence. Then 
$NG_{\kappa,\pi}^\mu \vert A_\kappa (\vec y) = f_{\vec y} (J)$ for some normal ideal $J$ on $P_\kappa (\lambda)$ extending $NG_{\kappa,\lambda}^\mu$.
\end{Obs}

{\bf Proof.} 
Set $K = NG_{\kappa,\pi}^\mu \vert A_\kappa (\vec y)$ and $J = p_\pi (K)$. Then by Observation 3.22, $J$ is normal, and moreover $K = f_{\vec y} (J)$. Finally by Fact 2.25 (ii), $NG_{\kappa,\lambda}^\mu \subseteq p_\pi (NG_{\kappa,\pi}^\mu) \subseteq p_\pi (K) = J$. 
\hfill$\square$ 

\begin{fact} {\rm (\cite{ForMag2})} Suppose that $\kappa = \nu^+$, where $\cf (\nu) \not= \mu$, $\lambda'$ is a cardinal greater than $\kappa$ of cofinality at least $\kappa$, and $K$ is a normal, fine ideal on $P_\kappa (\lambda')$ with the property that

\centerline{$\{x \in P_\kappa (\lambda') : \cf (\sup (x \cap \kappa)) = \cf (\sup x) = \mu\} \in K^+$.}

Then $K$ is not $(\lambda')^+$-saturated. 
\end{fact} 

\begin{Pro} Suppose that $\kappa = \nu^+$, where $\cf (\nu) \not= \mu$, and ${\mathcal B}_{\kappa,\lambda} (\kappa,\lambda^+)$ holds. Then there is an $(NG^\mu_{\kappa,\lambda}, I_{\kappa, \lambda})$-tower of length $\lambda^{++}$. 
\end{Pro} 

{\bf Proof.} Select a ${\mathcal B}_{\kappa,\lambda} (\kappa,\lambda^+)$-sequence $\vec y$, and set $K = NG_{\kappa,\pi}^\mu \vert A_\kappa (\vec y)$ and $J = p_\pi (K)$. By Facts 2.29 and 5.10, $K$ is not $(\lambda^+)^+$-saturated, so by Observation 4.20, there is a $(K, I_{\kappa, \lambda^+})$-tower of length $\lambda^{++}$. Now by (the proof of) Observation 5.9, $f_{\vec y} (NG^\mu_{\kappa,\lambda}) \subseteq f_{\vec y} (J) = K$. Moreover by Observation 3.14, $I_{\kappa, \lambda^+} \subseteq I_{\kappa, \lambda^+} \vert ran (f_{\vec y}) \subseteq f_{\vec y} (I_{\kappa, \lambda})$. By Observation 4.33, the desired conclusion follows.
\hfill$\square$ 

\medskip

To find towers in the other cases (i.e. when $\cf (\lambda) = \mu$ and $\kappa$ is either weakly inaccessible, or the successor of a cardinal of cofinality $\mu$), we proceed otherwise. For the same reasons as in the case when $\kappa \leq \cf (\lambda) < \lambda$, we will assume that $\mu = \omega$.

\medskip

\begin{fact} {\rm (\cite{Game})} Suppose that  $\cov (\lambda, \kappa, \kappa, \omega_1) = \lambda$. Then $NG^\omega_{\kappa,\lambda} = NS_{\kappa,\lambda} \vert T$ for some $T$.
\end{fact} 

\begin{Def}  An ideal $I$ on $\omega$ is a {\it P-point} if for any $F : \omega \rightarrow I$, there is $G \in I^\ast$ such that $G \cap F (j)$ is finite for all $j < \omega$. 
\end{Def}

\medskip

The following improves an earlier result of the author \cite{Nonsat1}.

\medskip

\begin{Pro} Suppose that
\begin{itemize}
\item $I$ is a P-point ideal on $\omega$. 
\item $\theta$ is a cardinal such that $\kappa < \theta \leq \lambda$ and $\cf (\theta) = \omega$.
\item $\langle \theta_i : i < \omega \rangle$ is a one-one sequence of regular infinite cardinals less than $\theta$ with supremum $\theta$ such that $\{i < \omega : \theta_i  \leq \xi\} \in I$ for all $\xi < \theta$.
\item $\lambda < \tcf(\prod_{i < \omega} \theta_i, <_I ) = \pi$. 
\item $\vec f = \langle f_\alpha: \alpha < \pi \rangle$  is an increasing, cofinal sequence in $(\prod_{i < \omega} \theta_i ,  <_I)$. 
\item $\psi : P_\kappa (\pi) \rightarrow \pi$ is defined by $\psi (x) =$ the least $\alpha$ such that $\{i :  \sup (x \cap \theta_i) \leq f_\alpha (i) \} \in I^\ast$.
\end{itemize}
Then letting $K$ be the $\pi$-complete ideal on $P_\kappa(\lambda)$ defined by 

\centerline{$K = \bigcup_{\alpha < \pi} P (\{ b \in P_\kappa (\lambda) : \psi (b) \leq \alpha \})$,}

the following hold :
\begin{enumerate}[\rm (i)]
\item There exists an $(NG^\omega_{\kappa,\lambda}, K)$-tower of length ${\frak b}_\pi$.
\item Suppose further that $\cov (\lambda, \kappa, \kappa, \omega_1) = \lambda$ and $cof (NS_{\kappa,\lambda}) \leq \pi $. Then for any $\eta < Depth ([\pi]^\pi, \searrow)$, there exists an $(NG^\omega_{\kappa,\lambda}, K)$-tower of length $\eta$.
\end{enumerate}
\end{Pro}

{\bf Proof.} 

(i) : By Proposition 4.31, we may find an $(NS_\pi \vert E^\pi_\omega, I_\pi)$-tower $\langle S_\alpha : \alpha < {\frak b}_\pi \rangle$. Set $X = \{ x \in P_\kappa (\pi) : \psi (x) = \sup x \}$ and for each $\alpha < \pi$, $A_\alpha = P_\kappa (\lambda) \cap \psi^{- 1} (S_\alpha)$. By a result of \cite{Menas}, $Sup (NG^\omega_{\kappa,\pi} \vert X) \subseteq NS_\pi \vert E^\pi_\omega$, so by Observations 4.2 (ii) and 4.15, $\langle Sup^{- 1} (S_\alpha) : \alpha < {\frak b}_\pi \rangle$ is an $(NG^\omega_{\kappa,\pi} \vert X, I_{\kappa, \pi})$-tower. We will show that $\langle A_\alpha : \alpha < {\frak b}_\pi \rangle$ is an $(NG^\omega_{\kappa,\lambda}, K)$-tower. Thus let $\alpha < \beta < \pi$.

\medskip

{\bf Claim 1. }  $p_\pi^{- 1} (X \cap (Sup^{- 1} (S_\alpha) \setminus Sup^{- 1} (S_\beta))) \subseteq A_\alpha \setminus A_\beta$. 

\medskip

{\bf Proof of Claim 1.}  Let $b \in p_\pi^{- 1} (X \cap (Sup^{- 1} (S_\alpha) \setminus Sup^{- 1} (S_\beta)))$. There must be $x \in X \cap (Sup^{- 1} (S_\alpha) \setminus Sup^{- 1} (S_\beta))$ such that $b = x \cap \lambda$. Then $\psi (b) = \psi (x) = \sup x$, so $\psi (b) \in S_\alpha \setminus S_\beta$, which completes the proof of the claim.

\medskip

Now $X \cap (Sup^{- 1} (S_\alpha) \setminus Sup^{- 1} (S_\beta)) \in (NG^\omega_{\kappa,\pi})^+$, so by Claim 1 and Fact 2.25 (ii), $A_\alpha \setminus A_\beta \in (NG^\omega_{\kappa,\lambda})^+$.

\medskip

{\bf Claim 2. }  $A_\beta \setminus A_\alpha \in K$.

\medskip

{\bf Proof of Claim 2.} Pick $\gamma < \pi$ with $S_\beta \setminus S_\alpha \subseteq \gamma$. Then $A_\beta \setminus A_\alpha = \{ b \in P_\kappa (\lambda) : \psi (b) \in S_\beta \setminus S_\alpha \} \subseteq \{ b \in P_\kappa (\lambda) : \psi (b) < \gamma \}$, which completes the proof of the claim, and that of (i).

\medskip

(ii) : By Fact 5.6, $cof (NG^\omega_{\kappa,\lambda}) \leq cof (NS_{\kappa,\lambda}) \leq \pi < {\frak b}_\pi$, so by a result of \cite{Menas}, there is $D \in NS_\pi^\ast$ such that $(\psi \vert P_\kappa (\lambda)) (NG^\omega_{\kappa,\lambda}) \subseteq I_\pi \vert (D \cap E_\omega^\pi)$. Now given $0 <\eta < Depth ([\pi]^\pi, \searrow)$, we may find an $(I_{\pi}, I_{\pi})$-tower $\langle S_\alpha : \alpha < \eta \rangle$ such that $S_\alpha \subseteq D \cap E_\omega^\pi$ for every $\alpha < \eta$. For $\alpha < \eta$, put $W_\alpha = P_\kappa (\lambda) \cap \psi^{- 1} (S_\alpha)$. Now let $\alpha < \beta < \eta$. Then $W_\alpha \setminus W_\beta$ equals $\{ b \in P_\kappa (\lambda) : \psi (b) \in S_\alpha \setminus S_\beta \}$, and lies consequently in $(NG^\omega_{\kappa,\lambda})^+$. On the other hand, $W_\beta \setminus W_\alpha$ equals $\{ b \in P_\kappa (\lambda) : \psi (b) \in S_\beta \setminus S_\alpha \}$, so it lies in $K$.
\hfill$\square$ 

\begin{Rmk} The P-pointness condition is annoying. There are however situations in which we do not have to worry about it. In fact, suppose that $\kappa < \theta \leq \lambda < \pi = \cf (\pi) < \min \{\pp^+(\theta), \theta^{+ \frak{p}} \}$, where $\cf (\theta) = \omega$ and $\frak{p}$ denotes the pseudo-intersection number. Then (see \cite{Shelah1}) for some countable set $A$ of uncountable regular cardinals with $\sup A = \theta$, $\tcf (\prod A / I) = \pi$, where $I$ is the noncofinal ideal on $A$.
\end{Rmk}

\begin{Cor} If $\cf (\lambda) = \omega$, then there exists an $(NG^\omega_{\kappa,\lambda}, I_{\kappa, \lambda})$-tower of length ${\frak b}_{\lambda^+}$.
\end{Cor}

{\bf Proof.} Use Fact 1.4.
\hfill$\square$ 


  \bigskip
\noindent Universit\'e de Caen - CNRS \\
Laboratoire de Math\'ematiques \\
BP 5186 \\
14032 Caen Cedex\\
France\\
Email :  pierre.matet@unicaen.fr\\


\end{document}